\documentclass[11pt,reqno]{amsart}
\def\squarebox#1{\hbox to #1{\hfill\vbox to #1{\vfill}}}

\newcommand{\Z}{{\mathbb Z}}

\newcommand{\R}{{\mathbb R}}
\newcommand{\C}{{\mathbb C}}
\newcommand{\N}{{\mathbb N}}

\theoremstyle{plain}

\newtheorem{thm}{Theorem}
  
\newtheorem{cor}{Corollary}
  
\newtheorem{lem}{Lemma}
  
\newtheorem{prop}{Proposition}

\newtheorem{rem}{Remark}

\begin{document}

\def\sn{{\bf S}^{n-1}}
\def\ts{\tilde{\sigma}}
\def\ss{{\mathcal S}}
\def\aa{{\mathcal A}}
\def\iu{\underline{i}}
\def\lc{{\mathcal L}}
\def\pxi{\phi + (\xi + i u)\psi_p}
\def\pp{\mathcal P}
\def\rr{\mathcal R}
\def\hU{\widehat{U}}
\def\mt{\Lambda}
\def\hLa{\widehat{\mt}}
\def\ep{\epsilon}
\def\tPhi{\widetilde{\Phi}}
\def\hR{\widehat{R}}
\def\cc{\mathcal C}
\def\oV{\overline{V}}
\def\uu{\mathcal U}
\def\tz{\tilde{z}}
\def\hz{\hat{z}}
\def\hd{\hat{\delta}}
\def\ty{\tilde{y}}
\def\hs{\hat{s}}
\def\hc{\hat{\cc}}
\def\hC{\widehat{C}}
\def\hh{\mathcal H}
\def\tf{\tilde{f}}
\def\of{\overline{f}}
\def\trr{\tilde{r}}
\def\tr{\tilde{r}}
\def\tts{\tilde{\sigma}}
\def\tVl{\widetilde{V}^{(\ell)}}
\def\tVj{\widetilde{V}^{(j)}}
\def\tVo{\widetilde{V}^{(1)}}
\def\tVj{\widetilde{V}^{(j)}}
\def\tPsi{\tilde{\Psi}}
 \def\tp{\tilde{p}}
 \def\tVjo{\widetilde{V}^{(j_0)}}
\def\tvj{\tilde{v}^{(j)}}
\def\tVjj{\widetilde{V}^{(j+1)}}
\def\tvl{\tilde{v}^{(\ell)}}
\def\tVt{\widetilde{V}^{(2)}}
\def\Lo{\; \stackrel{\circ}{L}}
\def\tg{\tilde{g}}

\def\ii{{\imath }}
\def\jj{{\jmath }}
\vspace*{0,8cm}
\def\saa{\Sigma_A^+}
\def\sa{\Sigma_A}
\def\sAA{\Sigma_{\aa}^+}
\def\Lip{\mbox{\rm Lip}}
\def\clip{C^{\mbox{\footnotesize \rm Lip}}}
\def\lip{\mbox{{\footnotesize\rm Lip}}}
\def\Vol{\mbox{\rm Vol}}
\def\lc{{\mathcal L}}

\def\ccm{\cc^{(m)}}
\def\oo{\mbox{\rm O}}
\def\ooo{\oo^{(1)}}
\def\oot{\oo^{(2)}}
\def\ooj{\oo^{(j)}}
\def\fo{f^{(1)}}
\def\ft{f^{(2)}}
\def\fj{f^{(j)}}
\def\wo{w^{(1)}}
\def\wt{w^{(2)}}
\def\wj{w^{(j)}}
\def\Vo{V^{(1)}}
\def\Vt{V^{(2)}}
\def\Vj{V^{(j)}}

\def\Uo{U^{(1)}}
\def\Ut{U^{(2)}}
\def\Ul{U^{(\ell)}}
\def\Uj{U^{(j)}}
\def\wl{w^{(\ell)}}
\def\Vl{V^{(\ell)}}
\def\Ujj{U^{(j+1)}}
\def\wjj{w^{(j+1)}}
\def\Vjj{V^{(j+1)}}
\def\Ujo{U^{(j_0)}}
\def\wjo{w^{(j_0)}}
\def\Vjo{V^{(j_0)}}
\def\vj{v^{(j)}}
\def\vl{v^{(\ell)}}

\def\gl{\gamma_\ell}
\def\id{\mbox{\rm id}}
\def\piU{\pi^{(U)}}
\def\bs{\bigskip}
\def\ms{\medskip}
\def\Int{\mbox{\rm Int}}
\def\diam{\mbox{\rm diam}}
\def\di{\displaystyle}
\def\dist{\mbox{\rm dist}}
\def\ff{{\mathcal F}}
\def\i{{\bf i}}
\def\pr{\mbox{\rm pr}}
\def\co{\; \stackrel{\circ}{C}}
\def\la{\langle}
\def\ra{\rangle}
\def\supp{\mbox{\rm supp}}
\def\Arg{\mbox{\rm Arg}}
\def\Int{\mbox{\rm Int}}
\def\II{{\mathcal I}}
\def\e{\emptyset}
\def\endofproof{{\rule{6pt}{6pt}}}
\def\con{\mbox{\rm const }}
\def\Box{\spadesuit}
\def\be{\begin{equation}}
\def\ee{\end{equation}}
\def\beqn{\begin{eqnarray*}}
\def\eeqn{\end{eqnarray*}}
\def\MM{{\mathcal M}}
\def\tmu {\tilde{\mu}}
\def\Pr{\mbox{\rm Pr}}
\def\Prf{\mbox{\footnotesize\rm Pr}}
\def\htau{\hat{\tau}}
\def\btau{\overline{\tau}}
\def\hr{\hat{r}}
\def\tF{\widetilde{F}}
\def\tG{\widetilde{G}}
\def\trho{\tilde{\rho}}

\def\Intu{\Int^u}
\def\Ints{\Int^s}
\def\hPsi{\hat{\Psi}}
\def\hp{\hat{p}}
\def\hJ{\hat{J}}

\def\W{{\mathcal W}}
\def\w{{\sf w}}
\def\sAA{\Sigma_{\aa}^+}
\def\sA{\Sigma_{\aa}}
\def\ff{{\mathcal F}}
\def\tm{\widetilde{m}}
\def\tpsi{\tilde{\psi}}
\def\tJ{\tilde{J}}

\def\tep{\tilde{\ep}}
\def\hw{\hat{w}}
\def\hGa{\widehat{\Gamma}}
\def\otau{\overline{\tau}}

\def\hm{h^{(m)}}
\def\hmo{h^{(m-1)}}
\def\ho{h^{(0)}}
\def\oh{\overline{h}}
\def\zm{z^{(m)}}
\def\zo{z^{(0)}}
\def\zmo{z^{(m-1)}}
\def\tpsi{\tilde{\psi}}
\def\opsi{\overline{\psi}}
\def\tphi{\tilde{\phi}}

\def\Lipe{\Lip_e}
\def\Lipf{\mbox{\rm \footnotesize Lip}}
\def\Lipfe{\mbox{\rm \footnotesize Lip}_e}




\title[Sharp large deviations]
{Sharp large deviations for some hyperbolic systems}
\author[V. Petkov]{Vesselin Petkov}
\address{Universit\'e Bordeaux I, Institut de Math\'ematiques de Bordeaux, 351,
Cours de la Lib\'eration, 
33405  Talence, France}
\email{petkov@math.u-bordeaux1.fr}
\author[L. Stoyanov]{Luchezar Stoyanov}
\address{University of Western Australia, School of Mathematics 
and Statistics, Perth, WA 6009, Australia}
\email{luchezar.stoyanov@uwa.edu.au}

\maketitle

\begin{abstract}
 We prove a sharp large deviation principle concerning intervals shrinking with 
sub-exponential speed for certain models involving the Poincar\'e map related to a Markov family
for an Axiom A flow restricted to a basic set $\Lambda$ satisfying some additional regularity assumptions.
\end{abstract}

\ms

\section{Introduction}
\renewcommand{\theequation}{\arabic{section}.\arabic{equation}}

\subsection{Main result}

Let $(X,f,\mu)$ be an ergodic dynamical system, where $f: X \to X$ is a diffeomorphism and $\mu$ an ergodic probability measure. 
For  an observable $\Psi: X \to \R$, Birkhoff's ergodic theorem says that
$$\frac{\Psi^n(x)}{n} = \frac{\Psi(x) + \Psi(f(x)) + \ldots + \Psi(f^{n-1}(x))}{n}$$
converges for almost all $x\in X$ with respect to $\mu$ to the mean value $M_{\Psi} = \int_X \Psi \, d\mu$ of $\Psi$
over $X$. So, if a closed interval $\Delta$ does not contain the mean $M_{\Psi}$, then the
measure of the set $\{ x\in X : \Psi^n(x)/n \in \Delta\}$ for $n$ sufficiently large should be small.
The theory of large deviations provides exponential bounds for such measures.  

For example it follows from general large deviation principles (see \cite{kn:Kif},
\cite{kn:Y}, \cite{kn:OP}) that if $X$ is a mixing basic set for an Axiom A diffeomorphism $f$,
$\Phi$ and $\Psi$ are H\"older continuous functions on $X$ with equilibrium states 
$m_\Phi$ and $m_\Psi$, respectively, and $m_\Psi$ is not the measure of maximal entropy of $f$
on $X$, then there exists a real-analytic {\it rate function} $J : \Int(\II_\Psi) \longrightarrow [0,\infty)$,
where
$$\II_\Psi = \left\{ \int \Psi\, d m : m \in \MM_X \right\} ,$$
$\MM_X$ is the set of all $f$-invariant Borel probability measures on $X$,  such that
\be
\lim_{\delta \to 0} \lim_{n \to \infty} \frac{1}{n} \log m_\Phi \left( \left\{ x\in X : \frac{\Psi^n(x)}{n}
\in (p - \delta , p + \delta)\right\}\right) = - J(p),\: \forall p \in \Int ({\mathcal I}_{\Psi}).
\ee
Since $m_\Psi$ is not the measure of maximal entropy, $\Psi$ is not cohomologous to a constant and the interval 
${\mathcal I}_{\Psi}$ is non trivial and $\Int({\mathcal I}_{\Psi}) \not= \emptyset.$ Moreover, $J(p) = 0$ if and only if $p = \int \Psi \, dm_\Phi$.

Many results on large deviations for hyperbolic (discrete and continuous) dynamical systems have been established in both
the uniformly hyperbolic case (see \cite{kn:Kif}, \cite{kn:Y}, \cite{kn:OP}, \cite{kn:L}, \cite{kn:W}, \cite{kn:G} and the 
references given there) and the non-uniformly hyperbolic case (\cite{kn:AP}, \cite{kn:RY}, \cite{kn:MN}).

For shrinking intervals $(p-\delta_n, p + \delta_n)$ with $\delta_n \to 0$ as $n \to \infty$ it follows that we have
an {\bf upper bound}
\be \label{eq:1.2}
 \lim\sup_{n \to \infty} \frac{1}{n} \log m_\Phi \left( \left\{ x\in X : \frac{\Psi^n(x)}{n}
\in (p - \delta_n , p + \delta_n)\right\}\right)\leq - J(p),\: \forall p \in \Int ({\mathcal I}_{\Psi}).
\ee
It is natural to study the question about the existence of a {\bf lower bound} in (\ref{eq:1.2}). 
Recently, Pollicott and Sharp (\cite{kn:PoS2}) obtained a result of this kind in the case of a
hyperbolic diffeomorphism $f : X \longrightarrow X$. Assuming that the H\"older continuous function
$\Psi$ satisfies a certain Diophantine condition related to three periodic orbits of $f$,  $m_\Psi$ is not the measure of maximal entropy of $f$,
and the sequence $\{\delta_n\}$ of positive numbers is such that $1/\delta_n = O(n^{1+\kappa})$ as $n \to \infty$
for some appropriately chosen $\kappa > 0$, they proved that
\be \label{eq:1.3}
\lim_{n \to \infty} \frac{1}{n} \log m_\Phi \left( \left\{ x\in X : \frac{\Psi^n(x)}{n}
\in (p - \delta_n , p + \delta_n)\right\}\right) = - J(p) 
\ee
for all $p \in \Int(\II_\Psi)$. As a consequence they derived a fluctuation theorem in a similar setup.

An apparently interesting question is whether one can go further and obtain the same lower bound as in (\ref{eq:1.3})
with a sequence $\{\delta_n\}$ converging must faster to $0$. Our aim in this article is to obtain a class of
examples where this holds in the case when $\delta_n \to 0$ with {\it sub-exponential speed}, i.e. when
\begin{equation} \label{eq:1.4}
\lim_{n \to \infty} \frac{\log \delta_n}{n} = 0\;.
\end{equation}
Moreover, we also show that for the class of functions we deal with, if 
$$\lim_{n \to \infty} \frac{\log \delta_n}{n} = -\alpha_0$$
for some sufficiently small $\alpha_0 > 0$, the asymptotic (\ref{eq:1.3}) is not true and we have a lower bound $- J(p) - \alpha_0.$ 
Thus our result in this situation is optimal and we indeed have sharp large deviations. 
To our best knowledge it seems that this is the first result with a precise limit different from $-J(p).$\\

Unlike \cite{kn:PoS2}, in our model the role of $X$ is played by the union of all rectangles in a Markov family
for an Axiom A flow restricted to a basic set $\Lambda$ and $f$ is just the corresponding Poincar\'e map. 

We now proceed to state our assumptions and main result precisely. Let $\varphi_t : M \longrightarrow M$ ($t\in \R$) be a $C^2$ 
flow on Riemannian manifold $M$ and let $\mt$ be a basic set for $\varphi_t$. It follows from the construction of Bowen \cite{kn:B} 
(cf. also  Ratner \cite{kn:Ra}) that there exists a Markov family $\rr = \{ R_i\}_{i=1}^k$ of rectangles 
$R_i = [U_i  , S_i ]$ of arbitrarily small size $\chi > 0$ for the restriction of the flow $\varphi_t$ to $\mt$ 
(see Section 2 for terminology and definitions).
Set $R = \cup_{i = 1}^k  R_i$. Let $\pp: R \longrightarrow R$ and $\tau : R \longrightarrow [0,\infty)$ be the corresponding
Poincar\'e map and {\it first return time}, respectively, so  that $\varphi_{\tau(x)}(x) = \pp(x)$.
We can then model $\varphi_t$ on $\mt$ by using the so called suspended flow on the suspension set
$R_\tau = \{ (x,t) : x\in R, 0 \leq t\leq \tau(x)\}$ (see e.g. Ch. 6 in \cite{kn:PP}). 

Let $F, G : \mt \longrightarrow \R$ be H\"older continuous functions. We will assume that the representative
of $G$ on $R_\tau$ is {\it constant on stable leaves}, i.e. on each set of the form
$\{ ([x,y], t) : y\in S_i\}$, where $i= 1, \ldots,k$, $x\in U_i$ and $t\in [0,\tau(x)]$.

Throughout this paper we assume the following 

\ms

\noindent
{\bf Standing Assumptions}: (A) {\it $\varphi_t$ is a mixing flow on a basic set $\mt$,
$\varphi_t$ and $\mt$ satisfy the conditions {\rm (LNIC)}, $(R_1)$
and $(R_2)$ stated in Sect. $2$ below and the local holonomy maps along stable laminations 
through $\mt$ are uniformly Lipschitz.}

\ms

(B) {\it  $\rr = \{ R_i\}_{i=1}^k$ is a fixed  Markov family of rectangles 
$R_i = [U_i  , S_i ]$ for the restriction of the flow $\varphi_t$ to $\mt$, chosen so that
the matrix $\aa = (a_{i, j})_{i,j = 1}^k$ related to $\rr$  is irreducible.}

\ms

(C) {\it $F : \mt \longrightarrow \R$ is a H\"older continuous function, while $G : \mt \longrightarrow \R$
is Lipschitz and its representative in  the suspension space $R_\tau$ is constant on stable leaves.}

\bs
For $T \geq 0$ and $x\in \mt$ set
$$G^T(x) =  \int_0^{T} G (\varphi_t(x))\, dt ,$$
and let
$$\II_0 =  \left\{ \int_R G^{\tau(x)} (x) \, d m(x) : m \in {\mathcal M_\pp} \right\} ,$$
where $\MM_\pp$ is the set of all $\pp$-invariant Borel probability measures on $R$.
For any function $h$ on $R$, $x\in R$ and an integer $n \geq 1$ we set
$$h^n(x) = h(x) + h(\pp(x)) + \ldots+ h(\pp^{n-1}(x)) .$$

Under the standing assumptions above,  in this paper we prove the following main result.

\begin{thm}  
There exists a  constant $\mu_0 > 0$ such that
for any Lipschitz  function $G > 0$ on $\mt$ with $\frac{\Lip(G)}{\min G} \leq \mu_0$, 
for which  $G^{\tau(x)}(x)$ is a non-lattice function on $R$, 
there exists a rate function $J : \II_0 \longrightarrow [0,\infty)$ with the following property:
for every H\"older continuous function $F$ on $\mt$ there exists a constant 
$\rho = \rho (F,G) \in (0,1)$ such that for any sequence $\{\delta_n\}$ of positive  numbers decreasing to zero with
\begin{equation} \label{eq:1.5}
\lim_{n \to \infty}\frac{\log \delta_n}{n} = -\alpha_0 \leq 0
\end{equation}
for some $0 \leq \alpha_0 \leq -\frac{\log \rho}{2}$, we have
\begin{equation} \label{eq:1.6}
\lim_{n\to \infty} \frac{1}{n} \log m \left( \left\{ x\in R : \frac{G^{\tau^n(x)}(x)}{n}
\in (p - \delta_n , p + \delta_n)\right\}\right) = - J(p) - \alpha_0,\:\forall p \in \Int(\II_0) ,
\end{equation}
where $m$ is the equilibrium state of the function $F^{\tau(x)}(x)$ on $R$.
In particular, for $\alpha_0 = 0$ we get 
\begin{equation} \label{eq:1.7}
\lim_{n\to \infty} \frac{1}{n} \log m \left( \left\{ x\in R : \frac{G^{\tau^n(x)}(x)}{n}
\in (p - \delta_n , p + \delta_n)\right\}\right) = - J(p) \: , \:\forall p \in \Int(\II_0) .
\end{equation}
\end{thm}

The rate function $J$ is explicitly defined in Sect. 3 below. The definition of `non-lattice' is provided
in Sect. 2.

It should be mentioned that in Theorem 1 it is enough to assume that $G$ is essentially Lipschitz
and $F$ is essential H\"older continuous (see the definitions in 1.2.2 below) -- the proof given below
works without any changes in that case. 

Notice that given any Lipschitz function $G$ on $\mt$, the function $\tG = G+c$ satisfies 
$\frac{\Lip(\tG)}{\min \tG} \leq \mu_0$ for any sufficiently large constant $c > 0$. Further remarks
on the assumptions in Theorem 1 are given in the next sub-section.

\ms

\subsection{On the range of applicability of the main result}

We begin with some remarks concerning the standing assumptions.


\subsubsection{\bf Remarks on the condition (A)}

(a) It is well-known that in general the maps $\mt \ni x \mapsto E^u(x)$ (or   $E^s(x)$) are only
H\"older continuous (see e.g \cite{kn:Ha} or \cite{kn:PSW}). The same applies to the so called local
stable and unstable holonomy maps (see Sect. 2 below for the definitions).
The following {\it pinching condition} implies stronger regularity properties of these maps. 

\ms

\noindent
{\sc (P)}:  {\it There exist  constants $C > 0$ and $0 < \alpha \leq \beta$ such that for every $x\in \mt$ we have
$$\frac{1}{C} \, e^{\alpha_x \,t}\, \|u\| \leq \| d\phi_{t}(x)\cdot u\| \leq C\, e^{\beta_x\,t}\, \|u\|
\quad, \quad  u\in E^u(x) \:\:, t > 0 \;,$$
for some constants $\alpha_x, \beta_x > 0$ depending on $x$ but independent of $u$ and $t$ with
$\alpha \leq \alpha_x \leq \beta_x \leq \beta$ and $2\alpha_x - \beta_x \geq \alpha$ for all $x\in \mt$.}

\ms

For example in the case of contact flows $\varphi_t$, it follows from the results 
in \cite{kn:Ha} that assuming (P), the map $\mt \ni x \mapsto E^u(x)$ is  $C^{1+\ep}$ with 
$\ep = 2 \alpha/\beta -1 > 0$ (in the  sense that this map has a linearization at any $x\in \mt$  that depends H\"older
continuously on $x$). The same applies to the map $\mt \ni x \mapsto E^s(x)$.

Notice that when $n = 2$ (then the local unstable manifolds are one-dimensional) the condition (P)
is always satisfied. Geodesic flows on manifolds $M$ of strictly negative curvature satisfy the pinching condition (P),
provided the sectional curvature is between $-K_0$ and $-K_0/4$ for some constants $K_0 > 0$ 
(\cite{kn:HP}).

For open billiard flows in the exterior of a compact set $K$ in $\R^n$ which is a finite
union of disjoint strictly convex domains satisfying a certain no eclipse condition,
the condition (P) is always satisfied when the minimal distance between  distinct connected 
components of $K$ is relatively large compared to the maximal sectional curvature of $\partial K$ (\cite{kn:St3}). 
In particular, for such billiards the standing assumption (A) is always satisfied (\cite{kn:St3}).

\ms

(b) As shown in \cite{kn:St2}, the non-integrability condition (LNIC) always holds for contact flows $\varphi_t$
when $\dim(M) = 3$, and also for transitive contact Anosov flows with Lispchitz stable/unstable holonomy maps.
In fact, in the latter case (LNIC), $(R_1)$ and $(R_2)$ always hold (\cite{kn:St5}).
In particular, for geodesic flows on compact locally symmetric spaces the standing assumption (A) is always
satisfied. 

\ms

(c) It is proved in \cite{kn:St4} that the conditions $(R_1)$ and $(R_2)$ hold under some rather general
assumptions for the flow $\phi_t$ and the basic set $\Lambda$. In particular this is always the case
under the pinching condition (P).
\ms

(d) We expect that a further progress in the analysis of the strong spectral estimates for the iterations of the Ruelle operator will
extend the setup of the dynamical systems for which we can apply our arguments. For example, it is natural to conjecture that 
for contact Anosov flow one should be able to obtain the result of Theorem 1 without the assumption (A).

\bs

\noindent
\subsubsection{\bf Remarks on the condition $\frac{\Lip(G)}{\min G} \leq \mu_0$ in Theorem 1} 



(a) Let $F, G : \longrightarrow \R$ be H\"older continuous functions, and let 
\be \label{eq:2.7}
\Phi(x) = \int_0^{\tau(x)} F (\varphi_t(x))\, dt \quad , \quad
\Psi(x) = \int_0^{\tau(x)} G (\varphi_t(x))\, dt\;
\ee
for all $x\in R$.  Then $\Phi$ and   $\Psi$ are {\it essentially H\"older} on $R$, i.e. for each $i \neq j$ 
they are H\"older continuous  on  $R_i \cap \pp^{-1}(R_j)$ (whenever this is non-empty). In fact, for the 
latter it is enough to assume that $F$ and $G$ are {\it essentially H\"older continuous}, i.e. they are 
H\"older continuous on each {\it block} 
$$B_{ij} = \{ \varphi_t(x) : x\in R_i \cap \pp^{-1}(R_j) \:, \: 0 \leq t \leq \tau(x) \} $$
with $R_i \cap \pp^{-1}(R_j) \neq \e$. In a similar way one defines {\it essentially continuous} and 
{\it essentially Lipschitz functions} on $\mt$ (with respect to the Markov family $\rr$).

It is easy to see that Theorem 1 can be stated in terms of the functions $\Phi$ and $\Psi$. Indeed, the measure
$m$ in Theorem 1 is $m_\Phi$, the equilibrium state of $\Phi$, and it is easy to show that 
$$G^{\tau^n(x)}(x) = \Psi^n(x)$$
for any $x\in R$ (see  subsection 2.1 below). Moreover, $\htau \leq \tau(x) \leq \htau_0$ ($x\in R$) for some constants $0 < \htau \leq \htau_0$,
so $\tau^n(x)/n \in [\htau, \htau_0]$ for all $x\in \R$. Finally, $\II_0 = \II_\Psi$, and the rate function $J$ in Theorem 1
is simply the rate function related to $\Psi$ (see (3.3) below). Thus, (1.6) is equivalent to
\begin{equation} \label{eq:1.9}
\lim_{n\to \infty} \frac{1}{n} \log m_\Phi \left( \left\{ x\in R : \frac{\Psi^n(x)}{n}
\in (p - \delta_n , p + \delta_n)\right\}\right) = - J(p) - \alpha_0 \: ,\:\forall p \in \Int(\II_\Psi) .
\end{equation}

(b) It is well-known that for any essentially continuous function $\Psi$ on $R$ there exists
an essentially continuous function $G$ on $\mt$ for which the second equality in (1.8) holds.
Indeed, fix a constant $A \in \R$ and a smooth function $\lambda : [0,1] \longrightarrow [0,1]$ such that
$\lambda(0) = 0$, $\lambda(1) = 1$ and $\lambda'(0) = \lambda'(1) = 0$, and set
$$G(\varphi_t(x)) = A + \left(\frac{\Psi(x)}{\tau(x)} -A\right) \, \lambda'(t/\tau(x)) \quad ,\quad
x\in R\, , \, 0 \leq t\leq \tau(x)\,.$$
One checks that $G$ satisfies (1.8) and $G(x) = G(\pp(x)) = A$ for all $x\in R$, which makes $G$
continuous on every block $B_{ij}$. Moreover, if $\tau$ is essentially Lipschitz, (which is the case under our assumptions; see Sect. 2)
$G$ is essentially Lipschitz (H\"older) whenever $\Psi$ is essentially Lipschitz (H\"older), and $G$ is constant
on stable leaves whenever $\Psi$ is.

\ms

(c)  Given an essentially Lipschitz function $\Psi$ on $R$, let $\Lipe (\Psi)$ be the smallest constant $L \geq 0$ so that
for any $x,y\in R_i \cap \pp^{-1}(R_j)$ we have $|\Psi(x) - \Psi(y)| \leq L d(x,y)$. (In a similar way we define $\Lipe(G)$
for an essentially Lipschitz function $G$ on $\mt$.)
Clearly, for any such $\Psi$ and any $\delta > 0$, adding a sufficiently large
constant $C > 0$ to $\Psi$ gives a function $\Psi + C$ satisfying 
$\frac{\Lipfe(\Psi+C)}{\min (\Psi +C)} \leq \delta$. Apart from this, the Ruelle operators
$L_{f -(\xi + \i b)\Psi}$ and $L_{f -(\xi + \i b)(\Psi+C)}$ are easily related, namely
$$\lc^n_{f -(\xi + \i b)(\Psi+C)}h = e^{- n(\xi + \i b) C}\, \lc^n_{f -(\xi + \i b)\Psi}\;.$$
However it is not clear whether one can find $G$ satisfying (1.8) with $\frac{\Lipfe(G)}{\min G} \leq \mu_0$,
assuming $\frac{\Lipfe(\Psi)}{\min \Psi } \leq \delta$ for some sufficiently small $\delta > 0$.

\ms


(d) It is easy to see that there is a non-trivial open set of essentially Lipschitz functions $\Psi$ on $R$ for which 
Theorem 1 applies. More precisely, for any constant $c > 0$ there is an open neighbourhood $V(c\tau)$ of $c\tau$ in 
the space $\clip_e(U)$ of essentially Lipschitz functions on $R$ constant on stable leaves
such that for any $\Psi\in V(c\tau)$  Theorem 1 applies.
Indeed, given $c > 0$, take $\delta = \delta(c) > 0$ small and let $V(c\tau)$ be the set of those
$\Psi = c\tau + \omega$, where $\Lipe(\omega) + \|\omega\|_{\infty} < \delta$. Given such $\Psi = c \tau + \omega$, define $G$ by
$$G(\varphi_t(x)) = c + \frac{\omega(x)}{\tau(x)} \, \lambda'(t/\tau(x)) \quad ,\quad
x\in R\, , \, 0 \leq t\leq \tau(x)\,,$$
where $\lambda$ is as in 1.2.2(b) above. Then (1.8) holds, and if $\delta$ is chosen sufficiently small, we have
$\min G \geq c/2$ and  $\Lipe(G) \leq C\delta$ for some constant $C > 0$ (depending on $\tau$ and $\lambda$ only),
so $\frac{\Lipfe(G)}{\min G} \leq \frac{2C\delta}{c}$, which can be made arbitrarily small choosing $\delta$ appropriately.

In comparison, one should remark that the Diophantine condition used in \cite{kn:PoS2} for the function $\Psi$ related 
to three periodic orbits is not an "open condition", although the functions satisfying this condition form a dense space.


 In  \cite{kn:W} Waddington proved an asymptotic of the form
\begin{equation}
m_{F}\left(\left\{ x \in \Lambda:\: G^T(x) - T q \in [\alpha, \beta]\right\}\right) \sim \Bigl(\int_{\alpha}^{\beta} e^{-\rho(q)t} dt\Bigr) 
\frac{C(q)}{\sqrt{2\pi \gamma''(\rho(p))}} \frac{e^{-T I(q)}}{\sqrt{T}},\: T \to + \infty ,
\end{equation}
where $m_F$ is the equilibrium state of $F$ on $\mt$ and $q  = \int_{\Lambda} G d m_{G + t F}.$ Here
$$-I(q) = \inf \{ \Pr(F + rG) - \Pr(F) - r q:\: r \in \R\}$$
and $\gamma(t) = \Pr(F  + tG) - \Pr(F),\: \gamma'(\rho(q)) = q$, where $\Pr = \Pr_\varphi$ is the topological pressure
with respect to the flow map $\varphi_1$ on $\mt$. This corresponds to large deviation for intervals $[\alpha/T, \beta/T]$ as $T \to \infty.$
Our Theorem 1 provides an asymptotic for
$$\frac{1}{n} \log m\left( \left\{ x \in R : G^{\tau^n(x)} - n p \in [-\delta'_n, \delta'_n]\right\} \right) \quad ,\quad  n \to \infty ,$$
where $\delta'_n =  n \delta_n \to 0$ sub-exponentially (or exponentially) fast, and $m = m_\Phi$ is the equilibrium
state of $\Phi(x) = F^{\tau(x)}(x)$ on $R$. It should be remarked that the ranges of Waddington's $q$'s and that of our $p$'s
are in general different, and so are the rate functions $I(q)$ in \cite{kn:W} and $J(p)$ in Theorem 1. It is natural to conjecture that an analogue of Theorem 1 holds for 
$$\frac{1}{n} \log \Bigl(m_F \{ x \in R:\: G^{\tau^n(x)} - \tau^n(x) q \in [-\delta_n,  \delta_n]\}\Bigr)$$
with $q = \int_{\Lambda} G d\mu$, $\mu$ being a $\varphi_t$-invariant probability measure, and $\delta_n$ satisfying (1.5). 
However this cannot be proved by using Theorem 2, since to apply the argument in Section 4 below, we need to know that the 
Ruelle operator ${\mathcal L}_{\Phi + (\xi + i u) (\Psi - q \tau)}$  is eventually contracting (see Theorem 2 below).  Since the function
 $\Psi(x) - q\tau(x)$ on $R$ is generated by $G(x) - q$, to use Theorem 2 we would need to have the condition
$\frac{{\rm Lip}_e\: (G)}{\min(G- q)} \leq \mu_0$ which could not be satisfied. Whether the latter condition can be omitted in the
assumptions of Theorem 2 is another open problem.


\section{Preliminaries}
\setcounter{equation}{0}

\subsection{Hyperbolic flows on basic sets}

Throughout this paper $M$ denotes a $C^2$ complete (not necessarily compact) 
Riemannian  manifold,  and $\varphi_t : M \longrightarrow M$ ($t\in \R$) a $C^2$ flow on $M$. A
$\varphi_t$-invariant closed subset $\mt$ of $M$ is called {\it hyperbolic} if $\mt$ contains
no fixed points  and there exist  constants $C > 0$ and $0 < \lambda < 1$ such that 
 there exists a $d\varphi_t$-invariant decomposition 
$T_xM = E^0(x) \oplus E^u(x) \oplus E^s(x)$ of $T_xM$ ($x \in \mt$) into a direct sum of non-zero linear subspaces,
where $E^0(x)$ is the one-dimensional subspace determined by the direction of the flow
at $x$, $\| d\varphi_t(u)\| \leq C\, \lambda^t\, \|u\|$ for all  $u\in E^s(x)$ and $t\geq 0$, and
$\| d\varphi_t(u)\| \leq C\, \lambda^{-t}\, \|u\|$ for all $u\in E^u(x)$ and  $t\leq 0$.
A non-empty compact $\varphi_t$-invariant hyperbolic subset $\mt$ of $M$ which is not a single 
closed orbit is called a {\it basic set} for $\varphi_t$ if $\varphi_t$ is transitive on $\mt$ 
and $\mt$ is locally maximal, i.e. there exists an open neighbourhood $V$ of
$\mt$ in $M$ such that $\mt = \cap_{t\in \R} \varphi_t(V)$. When $M$ is compact and $M$ itself
is a basic set, $\varphi_t$ is called an {\it Anosov flow}.

For $x\in \Lambda$ and a sufficiently small $\epsilon > 0$ let 
$$W^s_\ep(x) = \{ y\in M : d (\varphi_t(x),\varphi_t(y)) \leq \epsilon \: \mbox{\rm for all }
\: t \geq 0 \; , \: d (\varphi_t(x),\varphi_t(y)) \to_{t\to \infty} 0\: \}\; ,$$
$$W^u_\ep(x) = \{ y\in M : d (\varphi_t(x),\varphi_t(y)) \leq \epsilon \: \mbox{\rm for all }
\: t \leq 0 \; , \: d (\varphi_t(x),\varphi_t(y)) \to_{t\to -\infty} 0\: \}$$
be the (strong) {\it stable} and {\it unstable manifolds} of size $\epsilon$. Then
$E^u(x) = T_x W^u_\ep(x)$ and $E^s(x) = T_x W^s_\ep(x)$. 
Given $\delta > 0$, set $E^u(x;\delta) = \{ u\in E^u(x) : \|u\| \leq \delta\}$;
$E^s(x;\delta)$ is defined similarly. For any $A\subset M$ and  $I \subset \R$  denote
$\varphi_I(A) = \{\; \varphi_t(y)\; : \; y\in A, t \in I \; \}.$

From now on we will assume that $\mt$ is a basic set for $\varphi_t$. 

It follows from the hyperbolicity of $\mt$  that if  $\epsilon_0 > 0$ is sufficiently small,
there exists $\ep_1 > 0$ such that if $x,y\in \mt$ and $d (x,y) < \ep_1$, 
then $W^s_{\ep_0}(x)$ and $\varphi_{[-\ep_0,\ep_0]}(W^u_{\ep_0}(y))$ intersect at exactly 
one point $[x,y ] \in \mt$  (cf. \cite{kn:KH}). That is, there exists a unique $t\in [-\ep_0, \ep_0]$ 
such that $\varphi_t([x,y]) \in W^u_{\ep_0}(y)$. Setting $\Delta(x,y) = t$, defines the so called 
{\it temporal distance function} (\cite{kn:KB}, \cite{kn:D}). For $x, y\in \mt$ with 
$d (x,y) < \ep_1$, define
$\pi_y(x) = [x,y] = W^s_{\ep}(x) \cap \varphi_{[-\ep_0,\ep_0]} (W^u_{\ep_0}(y))\;.$
Thus, for a fixed $y \in \mt$, $\pi_y : W \longrightarrow \varphi_{[-\ep_0,\ep_0]} (W^u_{\ep_0}(y))$ 
is the {\it projection} along local stable manifolds defined on a small open neighbourhood $W$ of 
$y$ in $\mt$. Choosing $\ep_1 \in (0,\ep_0)$ sufficiently small, 
the restriction $\pi_y: \varphi_{[-\ep_1,\ep_1]} (W^u_{\ep_1}(x)) 
\longrightarrow \varphi_{[-\ep_0,\ep_0]} (W^u_{\ep_0}(y))$
is called a {\it local stable holonomy map\footnote{In a similar way one can define
holonomy maps between any two sufficiently close local transversals to stable laminations.}.} 
Combining such a map with a shift along the flow we get another local stable  holonomy  map
$\hh_x^y : W^u_{\ep_1}(x) \cap \mt  \longrightarrow W^u_{\ep_0}(y) \cap \mt$.
In a similar way one defines local holonomy maps along unstable laminations. It is well-known 
that the local holonomy maps are uniformly H\"older. Below we will assume that they are 
actually uniformly Lipschitz, i.e. choosing $\ep_1 > 0$ sufficiently small, there exists a 
constant $C_0 \geq 1$ such that
\be
d(\hh_x^y(u), \hh_x^y(v)) \leq C_0\, d(u,v)
\ee
for all $x, y \in \mt$ with $d(x,y) \leq \ep_1$, and all $u,v \in W^u_{\ep_1}(x)\cap \mt$.

Given $A\subset \mt$, denote by  $\diam(A)$ the {\it diameter} of $A$. We will say that
$A$ is an {\it admissible subset} of $W^u_{\ep}(z) \cap \mt$ ($z\in \mt$) if $A$ coincides with the closure
of its interior in $W^u_\ep(z) \cap \mt$. Admissible subsets of $W^s_\ep(z) \cap \mt$ are defined similarly.
Following \cite{kn:R} and \cite{kn:D}, a subset $R$ of $\mt$ will be called a {\it rectangle} if it has the form
$R = [U,S] = \{ [x,y] : x\in U, y\in S\}$, where $U$ and $S$ are admissible 
subsets of $W^u_\ep(z) \cap \mt$ and $W^s_\ep(z) \cap \mt$, respectively, for some $z\in \mt$. 
In what follows we will denote by 
$\Intu(U)$ the {\it interior} of $U$ in the set $W^u_\ep(z) \cap \mt$. In a similar way we 
define $\Ints(S)$, and then set $\Int(R) = [\Intu(U), \Ints(S)]$.
Given $\xi = [x,y] \in R$, set $W^u_R(\xi) = [U,y] = \{ [x',y] : x'\in U\}$ and
$W^s_R(\xi) = [x,S] = \{[x,y'] : y'\in S\} \subset W^s_{\ep_0}(x)$. The {\it interiors}
of these sets in the corresponding leaves are defined by $\Intu(W^u_R(\xi)) = [\Intu(U),y]$
and $\Ints(W^s_R(\xi)) = [x,\Ints(S)]$.

Let $\rr = \{ R_i\}_{i=1}^k$ be a family of rectangles with $R_i = [U_i  , S_i ]$,
$U_i \subset W^u_\ep(z_i) \cap \mt$ and $S_i \subset W^s_\ep(z_i)\cap \mt$, respectively, for some 
$z_i\in \mt$.   Set $R =  \cup_{i=1}^k R_i\; .$
The family $\rr$ is called {\it complete} if  there exists $T > 0$ such that for every $x \in \mt$,
$\varphi_{t}(x) \in R$ for some  $t \in (0,T]$.  The {\it Poincar\'e map} $\pp: R \longrightarrow R$
related to a complete family $\rr$ is defined by $\pp(x) = \varphi_{\tau(x)}(x) \in R$, where
$\tau(x) > 0$ is the smallest positive time with $\varphi_{\tau(x)}(x) \in R$.
The function $\tau$  is called the {\it first return time}  associated with $\rr$. 
A complete family $\rr = \{ R_i\}_{i=1}^k$ of rectangles in $\mt$ is called a 
{\it Markov family} of size $\chi > 0$ for the  flow $\varphi_t$ if $\diam(R_i) < \chi$ for all $i$ and: 
(a)  for any $i\neq j$ and any $x\in \Int (R_i) \cap \pp^{-1}(\Int (R_j))$ we have   
$\pp(\Int (W_{R_i}^s(x)) ) \subset \Int (W_{R_j}^s(\pp(x)))$ and 
$\pp(\Int (W_{R_i}^u(x))) \supset \Int (W_{R_j}^u(\pp(x)))$;
(b) for any $i\neq j$ at least one of the sets $R_i \cap \varphi_{[0,\chi]}(R_j)$ and
$R_j \cap \varphi_{[0,\chi]}(R_i)$ is empty.

The existence of a Markov family $\rr$ of an arbitrarily small size $\chi > 0$ for $\phi_t$
follows from the construction of Bowen \cite{kn:B} (cf. also  Ratner \cite{kn:Ra}). 

From now on we will assume that $\rr = \{ R_i\}_{i=1}^k$ is a fixed Markov family for  $\varphi_t$
of size $\chi < \ep_0/2 < 1$. Set  $U = \cup_{i=1}^k U_i\;.$
The {\it shift map} $\sigma : U   \longrightarrow U$ is given by
$\sigma  = \piU \circ \pp$, where $\piU : R \longrightarrow U$ is the {\it projection} along 
stable leaves. 
Notice that  $\tau$ is constant on each stable leaf $W_{R_i}^s(x) = W^s_{\ep_0}(x) \cap R_i$. 
For any integer $m \geq 1$
and any function $h : U \longrightarrow \C$ define $h^m : U \longrightarrow \C$ by
$h^m(u) = h(u) + h(\sigma(u)) + \ldots + h(\sigma^{m-1}(u))$. Then $\tau^m(x) = \tau(x) + \tau(\pp(x)) + \ldots + \tau(\pp^{m-1}(x)).$

Denote by $\widehat{R}$ the {\it core} of  $R$, i.e. the set of those $x\in R$ such that 
$\pp^m(x) \in \Int_{\mt}(R)$ 
for all $m \in \Z$. It is well-known (see \cite{kn:B}) that $\hR$ is a residual subset of $R$ and has full
measure with respect to any Gibbs measure on $R$. The set $\hU = U \cap \hR$ has similar properties.
Clearly in general $\tau$ is not continuous on $U$, however, under the standing assumption (A),
$\tau$ is {\it essentially Lipschitz} on $U$ in the sense that there exists a constant $L > 0$ such that if 
$x,y \in U_i \cap \sigma^{-1}(U_j)$ for some $i,j$, then $|\tau(x) - \tau(y)| \leq L\, d(x,y)$.
The same applies to $\sigma : \hU \longrightarrow \hU$.  Throughout we will mainly 
work with the restrictions of $\tau$ and $\sigma$ to $\hU$. Set $\hU_i = U_i \cap \hU$.

Given $z \in \mt$, let $\exp^u_z : E^u(z;\ep_0) \longrightarrow W^u_{\ep_0}(z)$  and
$\exp^s_z : E^s(z;\ep_0) \longrightarrow W^s_{\ep_0}(z)$ be the corresponding {\it exponential maps}.
A  vector $b\in E^u(z)\setminus \{ 0\}$ will be called  {\it tangent to $\mt$} at
$z$ if there exist infinite sequences $\{ v^{(m)}\} \subset  E^u(z)$ and 
$\{ t_m\}\subset \R\setminus \{0\}$
such that $\exp^u_z(t_m\, v^{(m)}) \in \mt \cap W^u_{\ep}(z)$ for all $m$, $v^{(m)} \to b$ and 
$t_m \to 0$ as $m \to \infty$. 
It is easy to see that a vector $b\in E^u(z)\setminus \{ 0\}$ is  tangent to $\mt$ at
$z$ if there exists a $C^1$ curve $z(t)$ ($0\leq t \leq a$) in $W^u_{\ep}(z)$ for some $a > 0$ 
with $z(0) = z$ and $\dot{z}(0) = b$ such that $z(t) \in \mt$ for arbitrarily small $t >0$.

The following  {\it local non-integrability condition} for $\varphi_t$ and $\mt$ was 
introduced in \cite{kn:St2}.

\bigskip

\noindent
{\sc (LNIC):}  {\it There exist $z_0\in \mt$,  $\ep_0 > 0$ and $\theta_0 > 0$ such that
for any  $\ep \in (0,\ep_0]$, any $\hz\in \mt \cap W^u_{\ep}(z_0)$  and any tangent vector 
$\eta \in E^u(\hz)$ to $\mt$ at $\hz$ with 
$\|\eta\| = 1$ there exist  $\tz \in \mt \cap W^u_{\ep}(\hz)$, 
$\ty_1, \ty_2 \in \mt \cap W^s_{\ep}(\tz)$ with $\ty_1 \neq \ty_2$,
$\delta = \delta(\tz,\ty_1, \ty_2) > 0$ and $\ep'= \ep'(\tz,\ty_1,\ty_2)  \in (0,\ep]$ such that
$$|\Delta( \exp^u_{z}(v), \pi_{\ty_1}(z)) -  \Delta( \exp^u_{z}(v), \pi_{\ty_2}(z))| \geq \delta\,  \|v\| $$
for all $z\in W^u_{\ep'}(\tz)\cap\mt$  and  $v\in E^u(z; \ep')$ with  $\exp^u_z(v) \in \mt$ and
$\la \frac{v}{\|v\|} , \eta_z\ra \geq \theta_0$,   where $\eta_z$ is the parallel 
translate of $\eta$ along the geodesic in $W^u_{\ep_0}(z_0)$ from $\hz$ to $z$. }

\bigskip

Set
$$B^u_T(x ,\ep) = \{ y \in W^u_\ep(x) \cap \mt : d(\varphi_t(x), \varphi_t(y)) \leq \ep \; , \; 0 \leq t \leq T\}\;.$$
Following \cite{kn:St2}, we will say that $\phi_t$ has a {\it regular distortion along unstable manifolds} over
the basic set $\mt$  if there exists a constant $\ep_0 > 0$ with the following properties:

\ms 
\noindent
($R_1$) For any  $0 < \delta \leq   \ep \leq \ep_0$ there exists a constant $R =  R (\delta , \ep) > 0$ such that 
$$\diam( \mt \cap B^u_T(z ,\ep))   \leq R \, \diam( \mt \cap B^u_T (z , \delta))$$
for any $z \in \mt$ and any $T > 0$.

\ms
\noindent
($R_2$) For any $\ep \in (0,\ep_0]$ and any $\rho \in (0,1)$ there exists $\delta  \in (0,\ep]$
such that for  any $z\in \mt$ and any $T > 0$ we have
$\diam ( \mt \cap B^u_T(z ,\delta))   \leq \rho \; \diam( \mt \cap B^u_T (z , \ep))\;.$

\bs

Let $B(\hU)$ be the {\it space of  bounded functions} $g : \hU \longrightarrow \C$ with its 
standard norm  $\|g\|_{\infty} = \sup_{x\in \hU} |g(x)|$. Given a function $g \in B(\hU)$, the  
{\it Ruelle transfer operator } $\lc_g : B(\hU) \longrightarrow B(\hU)$ is defined by 
$\di (\lc_gh)(u) = \sum_{\sigma(v) = u} e^{g(v)} h(v)\;.$
If $g \in B(\hU)$ is Lipschitz on $\hU$, then  $\lc_g$ preserves the space $\clip(\hU)$
of {\it Lipschitz functions} $g: \hU \longrightarrow \C$.
Given $\alpha > 0$, let $C^\alpha (U)$ be the space of H\"older continuous
functions $G: U \longrightarrow \C$ with
$$|h|_\alpha = \sup \left\{ \frac{|h(x) -h(y)|}{(d(x,y))^\alpha} : x,y \in U_i \;,\;
x\neq y\;,\; i = 1, \ldots,k\right\} < \infty\;.$$
We will consider $C^\alpha(U)$ with the norm
$\| h\|_{\alpha} = \|h\|_\infty + |h|_\alpha$.

The hyperbolicity of the flow on $\mt$ and the standing assumption (A) imply the existence of
 constants $c_0 \in (0,1]$ and $\gamma_1 > \gamma > 1$ such that
\begin{equation}
c_0 \gamma^m\; d (u_1,u_2) \leq 
d (\sigma^m(u_1), \sigma^m(u_2)) \leq \frac{\gamma_1^m}{c_0} d (u_1,u_2)
\end{equation}
whenever $\sigma^j(u_1)$ and $\sigma^j(u_2)$ belong to the same  $U_{i_j}$ 
for all $j = 0,1 \ldots,m$.  

\subsection{Non-lattice condition}

First  it easy to see that $G^{\tau^n(x)} = \Psi^n(x).$ Indeed,
$$\int_0^{\sum_{j=-1}^{n-1} \tau(\pp^j(x))} G(\varphi_t(x)) dt 
= \sum_{k = -1}^{n-2}\int_{\sum_{j=-1}^{k}\tau(\pp^j (x))}^{\sum_{j=-1}^{k+1} \tau(\pp^j(x))} G(\varphi_t(x)) dt,$$
where by convention $\tau(\pp^{-1}(x)) = 0$. On the other hand, setting $t = \sum_{j= -1}^{k}\tau(\pp^j (x))+ s$, we get
$$\int_{\sum_{j=-1}^{k}\tau(\pp^j (x))}^{\sum_{j=-1}^{k+1} \tau(\pp^j(x))} G(\varphi_t(x)) dt 
= \int_0^{\tau(\pp^{k+1}(x))} G(\varphi_{\sum_{j=-1}^{k}\tau(\pp^j (x))+ s}(x)) ds $$
$$= \int_0^{\tau(\pp^{k+1}(x))} G(\varphi_s (\pp^{k + 1}(x)) ds $$
and
$$\int_0^{\sum_{j=-1}^{n-1} \tau(\pp^j(x))} G(\varphi_t(x)) dt 
= \sum_{k=-1}^{n-2} \int_0^{\tau(\pp^{k+1}(x))} G(\varphi_s(\pp^{k+1}(x)) ds$$
$$= \Psi(x) + ...+ \Psi(\pp^{n-1}(x)) = \Psi(x) + \ldots \Psi(\sigma^{n-1}(x)) = \Psi^n(x) .$$

A H\"older continuous function $g(x)$ on $R$ is called {\it non-lattice} if there do not exist constant
$a$, a H\"older continuous function $h$ on $R$ and a bounded integer-valued function $Z$ on $R$ so that
$$g(x) = (h \circ \pp)(x) - h(x)  + a +  Z(x) \quad ,\quad  x \in R . $$
Notice that such a function $Z$ can only have a finite range.
If the function $G^{\tau(x)}$ is lattice with some constant $a$, then for every periodic orbit $\gamma$ of $\varphi_t$ issued from 
$x \in R$ with $\pp^n (x) = x$   we have
\begin{equation}
G^{\tau^n(x)} - an = \Psi^n (x) - an = \int_{\gamma} G(\varphi_t) dt  - an \in \Z.
\end{equation}
The above condition on $\Psi$ is the same as the lattice condition introduced in \cite{kn:PoS2}.
We show below that the condition (2.3) is related to another one. 
\ms

Recall that two (essentially) continuous functions $F$ and $G$ on $\mt$ are called {\it cohomologuous} ($F\sim G$)
if there exists an (essentially) continuously differentiable function $H$ on 
$\mt$ with $F - G = H'$ (see \cite{kn:La} or \cite{kn:W}). $H$ is called (essentially) {\it continuously differentiable} if there exists an
(essentially) continuous function $H'$ on $\mt$ such that
$$H'(x) = \lim_{t\to 0} \frac{H(\varphi_t(x)) - H(x)}{t} .$$
Notice that if $\Phi$ and $\Psi$ are defined by (1.8) and $F \sim G$, then $F - G = H'$ for some $H'$ as above, so
$$\Phi(x) - \Psi(x) = \int_0^{\tau(x)} H'(\varphi_t(x)) \, dt = \int_0^{\tau(x)} \frac{d}{dt} H(\varphi_t(x)) \, dt
= H(\varphi_{\tau(x)}(x)) - H(x) = H(\pp(x)) - H(x)$$
for all $x\in R$, and therefore $\Phi \sim \Psi$. Conversely (see \cite{kn:La}) if $\Phi \sim \Psi$, then $F \sim G$.
We can now express the non-lattice condition on $\Psi$ in terms of the function $G$. Assume for a moment that
$\Psi$ is lattice on $R$, i.e. there exist a constant $a\in \R$ and a bounded integer-valued 
function $Z$ on $R$ so that $\Psi \sim a + Z$. 
Let $G_0$ be a fixed essentially Lipschitz function on $\mt$ such that $G_0^{\tau(x)}(x) = 1$ for all $x\in R$, e.g. define 
$$G_0(\varphi_t(x)) = \frac{1}{\tau(x)}  \, \lambda'(t/\tau(x)) \quad ,\quad
x\in R\, , \, 0 \leq t\leq \tau(x)\,,$$
where $\lambda$ is as in 1.2.2(b). Setting 
\begin{equation} \label{eq:2.4}
M(\varphi_t(x)) = Z(x),\: \forall x\in R,\:  0 \leq t < \tau(x),
\end{equation}
 we get a
bounded integer-valued function on $\mt$, and $a + Z(x) = (aG_0 + M G_0)^{\tau(x)}(x)$ for all $x\in R$. Thus,
$\Psi \sim a + Z$ means that $G \sim aG_0 + MG_0 = (a+M) G_0$. Therefore a {\bf sufficient condition} for $\Psi$ to be
non-lattice is that there do not exist a constant $a\in \R$ and a bounded integer-valued function $M$ on $\mt$ satisfying 
(\ref{eq:2.4}) so that $G \sim (a+M) G_0$. This can be used to define a non-lattice condition for $G$, although it is a bit stronger than what is
necessary to make it equivalent to $\Psi$ being non-lattice on $R$.

\section{Eventually contracting Ruelle operators}
\setcounter{equation}{0}

\subsection{Topological pressure and rate functions}

Let $\varphi_t$ be a $C^2$ flow on a Riemann manifold $M$, $\mt$ be a mixing basic set for $\phi_t$,
and $\rr = \{ R_i\}_{i=1}^k$ be a fixed Markov family for  $\varphi_t$ of size $\chi < \ep_0/2 < 1$.

Denote by $\Pr_\pp(h)$ the {\it topological pressure} of a continuous function $h$ on  $R$ with respect to the map $\pp$ on $R$). 
We will often write just $\Pr(h)$ instead of $\Pr_\pp(h)$.
Recall that for any continuous function $\Psi : R \longrightarrow \R$,
\be
\Pr_{\pp}(\Psi) = \sup_{m \in {\mathcal M_{\pp}}} \Bigl( h_{\pp}(m) + \int_R \Psi d m \Bigr) ,
\ee
where ${\mathcal M}_{\pp}$ is the set of all $\pp$-invariant probability measures on $R$ 
and $h_{\pp}(m)$ is the measure theoretic entropy of $m$ with respect to $\pp$. 

Given a continuous function $H$ on $\mt$, set
$$H^\tau(x) = \int_0^{\tau(x)} H(\varphi_t(x))\, dt \quad , \quad x\in R .$$

Let $m$ be a $\pp$-invariant probability measure on $R$. 
There exists a unique $\varphi_t$-invariant probability measure $\mu$ on $\mt$ such that
\be
\int_\mt H\, d\mu = \frac{\int_R \left( \int_0^{\tau(x)} H(\varphi_t(x)\, dt\right)\, dm(x)} {\int_R \tau \, dm}
= \frac{\int_R H^{\tau} (x)\, dm(x)} {\int_R \tau \, dm}
\ee
for any continuous function $H$ on $\mt$ (see e.g. Ch. 6 in \cite{kn:PP}). 
The map $m \mapsto \Omega(m) = \mu$ is onto. Moreover, we have $\Pr_\pp(H^\tau - \Pr_\varphi(H)\,\tau) = 0$.
If $m = m_{H^\tau - \Prf_\varphi(H)\,\tau}$ is the {\it equilibrium state} of 
$H^\tau - \Pr_\varphi(H)\,\tau$ on $R$, then $\mu = \Omega(m)$ coincides with the equilibrium state $\mu_H$ 
of the function $H$ on $\mt$ (see Ch. 6 in \cite{kn:PP}).

For a H\"older continuous function $\Psi$ on $R$ set
$$
\II_\Psi = \left\{ \int_R \Psi\, d m : m \in {\mathcal M}_{\pp} \right\} .$$


It follows from the Large Deviation Theorem in \cite{kn:Kif} that if $m_\Psi$ is not 
the measure of maximal entropy for $\pp$, then there exists a real analytic function
$J : \Int(\II_\Psi) \longrightarrow [0,\infty)$ such that $J(p) = 0$ iff
$p = \int_R \Psi\, d m_\Phi$ for which (1.2) holds. More precisely, we have
\be
-J(p) = \inf \{\Pr_\pp(\Phi + q \Psi) - \Pr_\pp(\Phi) - qp: \: q \in \R\} .
\ee
It is also known that
\be
\left[\frac{d }{dq}\, \Pr_{\pp}(\Phi + q\Psi)\right]_{q= \eta} = \int_R \Psi \, dm_{\Phi+\eta \Psi} .
\ee
Moreover, for any $p \in  \II_\Psi$ there exists a unique number $\xi = \xi_p$ such that
\be
-J(p) = \Pr_{\pp} (\Phi + \xi\, \Psi) - \Pr_\pp(\Phi) -\xi \, p \quad , \quad p = \int_R \Psi\, d m_{\Phi + \xi \, \Psi} .
\ee

Let $F,G : \longrightarrow \R$ be H\"older continuous functions.
Under the standing assumptions,  define the functions $\Phi, \Psi : R \longrightarrow \R$ by
\be \label{eq:1.5}
\Phi(x) = \int_0^{\tau(x)} F (\varphi_t(x))\, dt \quad , \quad
\Psi(x) = \int_0^{\tau(x)} G (\varphi_t(x))\, dt .
\ee
Then  $\Psi$ is {\it essentially Lipschitz on $R$}, i.e.  for each $i \neq j$ it is Lipschitz on  $R_i \cap \pp^{-1}(R_j)$. Similarly, 
$\Phi$ is {\it essentially H\"older continuous} on $R$. This implies (see e.g. \cite{kn:PP}) that
there is a well-defined  {\it equilibrium state} $m_{\Phi}$ of the function $\Phi$ on $R$,
while for $\Psi$ the consequence is that the Ruelle transfer operators of the form
${\mathcal L}_{f+ s \Psi}$ with $s\in \C$ and $f$ H\"older continuous on $R$ and constant on stable leaves,
are well-defined in appropriately chosen spaces of H\"older continuous functions
(see Sect. 2  for the main definitions). 

As mentioned in Sect. 1.2, (1.6) in Theorem 1 is  equivalent to (1.9). So, this is what we are going
to prove in the next section.

\ms

\noindent
{\bf Remark 1.} Notice that for the proof of (1.9) it is enough to consider the case 
\be
\Pr_\pp(\Phi) = 0 ,
\ee
since we can change $\Phi$ by adding a constant and this preserves the measure $m_{\Phi}$.
Moreover, it is enough to prove (1.9) in the case when $\Phi$ is constant on stable leaves in $R$, and
\be
\int_R \Psi \, d\mu_\Phi = 0\;.
\ee
Indeed,  it follows from Sinai's Lemma (c.f. e.g. Ch. 1 in \cite{kn:PP}) that $\Phi$
is cohomologues to a H\"older continuous function $\tPhi$ on $R$ which is constant on stable leaves.
Then $m_{\Phi} = m_{\tPhi}$, so we can replace $\Phi$ by $\tPhi$. For the second statement above,
given functions $\Phi$ and $\Psi$ as in Theorem 1, set $\hPsi = \Psi - c$, where
$c = \int_R \Psi \, d\mu_\Phi$. Then $\int_R \hPsi\, d\mu_\Phi = 0$. Recall  that the function
$J : \Int(\II_\Psi) \longrightarrow \R$ is given by (3.3).

Let $\hJ : \Int(\II_{\hPsi}) \longrightarrow [0,\infty)$
be the corresponding functional for $\hPsi$, i.e.
$$- \hJ(\hp) = \inf \{\Pr_\pp(\Phi + q \hPsi) - q \hp : \: q \in \R\} \quad , \quad \hp \in \II_{\hPsi}\;.$$
One checks immediately that $\II_{\hPsi} = \II_{\Psi} - c$ and for any $\hp = p - c \in \II_{\hPsi}$,
using the properties of pressure,  for any $q\in \R$ we get
$$\Pr_\pp(\Phi + q \hPsi) - q \hp
= \Pr_\pp(\Phi + q (\Psi-c)) - q (p-c) = \Pr_\pp(\Phi + q\Psi) -qp\;.$$
Thus, $\hJ(\hp) = J(p)$. Hence if (1.9) holds with $\Psi$, $J$ and $p \in \II_\Psi$ replaced by
$\hPsi$, $\hJ$ and $\hp \in \II_{\hPsi}$, then (1.9) holds in its present form, as well.

\subsection{Ruelle transfer operator}

Given a Lipschitz real-valued function $f$  on $\hU$, set $\tf  = f - P\Psi$, where $P = P_f\in \R$ is the unique 
number such that the topological pressure $\Pr_\sigma(\tf)$ of $\tf$ with respect to $\sigma$ is 
zero (cf. e.g. \cite{kn:PP}). For $a, b\in \R$, one defines the {\it Ruelle transfer operator}
$$\lc_{\tf-(a+\i b)\Psi} : \clip (\hU) \longrightarrow \clip (\hU)$$ 
in the usual way (cf. Sect. 2 above).

We will say  that the {\it Ruelle transfer operators related to $\Psi$ and the function $f$ on $\hU$ are eventually contracting} 
if for every $\epsilon > 0$ there exist constants $0 < \rho < 1$, $a_0 > 0$ and  $C > 0$ such that if $a,b\in \R$  
satisfy $|a| \leq a_0$ and $|b| \geq 1/a_0$, then for every integer $m > 0$ and every  $h\in \clip (\hU)$ we have
$$\|\lc_{f -(P_f+a+ \i b)\Psi}^m h \|_{\lip,b} \leq C \;\rho^m \;|b|^{\ep}\; \| h\|_{\lip,b}\; ,$$
where the norm $\|.\|_{\lip,b}$ on $\clip (\hU)$ is defined by 
$\| h\|_{\lip,b} = \|h\|_{\infty} + \frac{\lip(h)}{|b|}$.
This implies in particular that the spectral radius  of $\lc_{f-(P_f+ a+\i b)\Psi}$ on $\clip(\hU)$ does not exceed  $\rho$.

The following theorem is one of the main ingredients in the proof of Theorem 1.

\begin{thm} Under the standing assumptions, let $\Psi : R \longrightarrow \R$ be defined by $(\ref{eq:1.5}).$
Then there exists a constant $\mu_0 > 0$ such that if 
$\frac{\Lipf(G)}{\min G} \leq \mu_0$, then for any H\"older continuous real-valued function 
$f$  on $\hU$ the Ruelle transfer operators related to  $\Psi$ and $f$  are eventually contracting.
\end{thm}

Theorem 1 is derived in Sect. 4, while Theorem 2 is proved in Sect. 5. using some ideas from 
\cite{kn:St2}.

\section{Proof of Theorem 1}
\setcounter{equation}{0}

We will use the notation and the assumptions in Sections 2 and 3. Let $\Phi$ and $\Psi$ be as in (3.6). 
As in the beginning of Sect. 3, we may assume (3.7) and (3.8). In this section we will always consider
pressure with respect to $\pp$, so for brevity we will write $\Pr$ instead of $\Pr_\pp$. Let $J$ be defined by (3.3).
Consider a sequence $\{\delta_n\}_{n \in \N},\:\delta_n > 0,\: \delta_n \to 0$, such that 
(\ref{eq:1.5}) holds and let $\epsilon_n = n \delta_n$.

Fix an arbitrary $p\in \Int(\II_\Psi)$ and set $\Psi_p = \Psi- p$.
As in \cite{kn:PoS2}, it is enough to prove a modified result concerning a sequence of the form
\be
\rho(n) = \int_{U} \chi_n(\Psi_p^n(x)) \; d \mu_{\Phi}\;,
\ee
where  $\chi \in C_0^k(\R: \R)$ is a {\it fixed  cut-off function} and
\be
\chi_n(x) = \chi(\epsilon_n^{-1} x) \quad , \quad x\in \R\;.
\ee

\begin{prop} Under the assumptions of Theorem $2$, we have
$$\lim_{n\to \infty} \frac{1}{n} \, \log \rho(n) = - J (p)- \alpha_0\;.$$
\end{prop}

Theorem 1 follows immediately from Proposition 1 as shown in \cite{kn:PoS2} choosing two functions 
$\chi_{-}(x), \: \chi_{+}(x) \in C_0^k(\R: \R)$ so that  $\chi_{-}(x) \leq {\bf 1}_{[-1, 1]}(x) \leq \chi_{+}(x)$ and 
$$2 - \eta \leq \int_R \chi_{-}(x)dx  \leq \int_R \chi_{+}(x) dx \leq 2 + \eta,\: \eta > 0.$$
The rest of this section is devoted to the proof of Proposition 1.

For the given functions $\Phi$ and $\Psi$, Theorem 2 implies the following.

\begin{cor} For any $\delta > 0$ and any $a  > 0$ there exist $\rho \in (0,1)$, 
$a_0 > 0$ and $A_\delta > 0$ (depending on $a$ as well) such that
\begin{equation} 
\|\lc_{\Phi + (\xi + iu)\Psi}^n 1\|_{\infty} \leq A_\delta\; \rho^n |u|^{\delta} e^{n\Prf (\Phi + \xi\, \Psi)}\;,
\end{equation}
for all integers $n \geq 1$, all $u, \xi \in \R$ with $|u|\geq 1/a_0$ and $|\xi| \leq a$.
\end{cor}

Let $\xi_p$ be the unique real number such that 
$$- J(p) = \inf_{q \in \R}\{\Pr(\Phi + q \Psi_p\} = \Pr(\Phi + \xi_p \Psi_p).$$
Therefore
$$\frac{d \Pr(\Phi + q \Psi_p)}{dq} \big\vert_{q = \xi_p} = \int \Psi_p \; d\mu_{\Phi + \xi_p \Psi_p} = 0\,.$$
Moreover, since $\Psi$ and $\Psi_p$ are non-lattice, we have
$$\frac{d^2 \Pr(\Phi + (\xi + iu) \Psi_p)}{du^2}\big\vert_{u = 0} = -\sigma^2 < 0.$$ 
The operator $\lc_{\Phi + \xi_p\Psi_p}$ has a simple eigenvalue $\lambda_{\xi} = e^{\Prf(\Phi + \xi_p \Psi_p)}$,
and so for all sufficiently small $u \in \C$ the operator $\lc_{\Phi + (\xi_p + i u)\Psi}$ has a simple eigenvalue
$e^{\Prf(\Phi + (\xi_p + iu)\Psi)}$ and the rest of the spectrum of $\lc_{\Phi + (\xi_p + iu)\Psi}$ is contained in 
a disk of radius $\theta \lambda_{\xi}$ with some $0 < \theta < 1.$

Clearly for the Fourier transform $\hat{\chi}$ of $\chi$ we get
$\hat{\chi}_n(u) = \epsilon_n \hat{\chi}(\epsilon_n u).$  
Then
$$\Bigl \{ x:\: \frac{\Psi^n(x)}{n} \in (p - \delta_n, p + \delta_n)\Bigr\} = \{ x:\: \Psi_p^n(x) \in (-\epsilon_n, \epsilon_n)\}.$$
Set $\xi = \xi_p$ and $\omega_n(y) = e^{-\xi y}\chi_n(y).$ We need the following lemma established in \cite{kn:PoS2}.

\begin{lem}  If $\Phi$ is normalized so that $\Pr(\Phi) = 0$, then
$$\int e^{(\xi + i u)\Psi_p^n(x)} d \mu_{\Phi}(x) = 
\int \lc^n_{\Phi + (\xi + iu)\Psi_p} 1(x) \; d \mu_{\Phi}(x)\;.$$
\end{lem}

\def\pxi{{\footnotesize{\Phi + (\xi + i\, u) \Psi_p}}}

Using the lemma and applying the Fourier transform, we have
$$\rho(n) = \frac{1}{2 \pi} \int_{-\infty}^{\infty} \Bigl( \int e^{iu \Psi_p^n(x)}  d \mu_{\Phi} (x)\Bigr) \hat{\chi}_n(u) du$$
$$ = \frac{1}{2\pi} \int_{-\infty}^{\infty} \Bigl( \int \lc_{\Phi + (\xi + i u)\Psi_p}^n 1 (x)  d\mu_{\Phi}(x)\Bigr) \hat{\omega}_n(u) du$$
$$= \frac{\epsilon_n}{2\pi} \int_{-\infty}^{\infty} \Bigl( \int \lc_{\Phi + (\xi + i u)\Psi_p}^n 1 (x) \;
d\mu_{\Phi}(x)\Bigr) \hat{\chi}(\epsilon_n (u - i \xi)) du.$$
We choose $a > 0$ sufficiently small and changing the coordinates on $(-a, a)$ to $v = v(u)$, we write
$$e^{\Prf(\pxi)} = \lambda_{\xi}(1 - v^2 + iQ(v)) \;,$$
where $Q(v)$ is real valued and $Q(v) = {\mathcal O}(|v|^3).$
The analysis in subsection 4.1 in \cite{kn:PoS2} yields
$$I_1 = \int_{-a}^a \Bigl( \int \lc_{\pxi}^n 1(x) d\mu_{\phi}(x)\Bigr) \hat{\omega}_n(u) du 
=  \frac{\hat{\chi}(0)}{\sigma} \frac{\epsilon_n}{\sqrt{2 \pi n}} \lambda_{\xi}^n 
+ {\mathcal O} \Bigl(\frac{\epsilon_n \lambda_{\xi}^n}{n}\Bigr).$$
Next, we consider the integral
$$I_2 = \epsilon_n\int_{a < |u| \leq c}\Bigl(\int \lc_{\pxi}^n 1(x) d\mu_{\Phi}(x)\Bigr) 
\hat{\chi}(\epsilon_n(u - i\xi)) du$$
with $c \gg 1$ sufficiently large. Since $\Psi_p$ is non lattice,  for $0 < a \leq |u| \leq c$
the operator ${\mathcal L}_{\Phi + (\xi + iu)\Psi_p}$ has no eigenvalues $\mu$ with  $|\mu| = \lambda_{\xi}$ 
(see for instance \cite{kn:PP}), and the spectral radius of $\lc_{\Phi + (\xi + i u)\Psi_p }$ 
is strictly less than $\lambda_{\xi}$. Thus, there exists $\beta = \beta(a, c),\: 0 < \beta < 1$, 
such that for $n \geq N(a, c)$ we have
\begin{equation} \label{eq:4.4}
\|{\mathcal L}^n_{\Phi + (\xi + i u) \Psi_p}\| \leq \beta^n\lambda_{\xi}^n.
\end{equation}
On the other hand,
\begin{equation} \label{eq:4.5}
|\hat{\chi}(\epsilon_n(u - i\xi))| \leq C_k \frac{e^{c_0 |\epsilon_n \xi|}}{\epsilon_n^k |u|^k}, \:|u| \geq a,\: \forall k \in \N\;,
\end{equation}
with $c_0 > 0$ depending on the support of $\chi.$ Applying (\ref{eq:4.4}) and (\ref{eq:4.5}) with 
$k = 0$, for large $n$ we get
$$I_2 = {\mathcal O} \Bigl(\frac{\epsilon_n \lambda_{\xi}^n}{n}\Bigr).$$

Now consider
$$I_3 = \epsilon_n \int_{|u| > c} \Bigl( \int \lc_{\pxi}^n 1(x) d\mu_{\Phi}(x)\Bigr)  \hat{\chi}(\epsilon_n(u - i\xi)) du\;.$$
We are going to use the spectral estimate (4.3).
Fix $0 < \delta \leq 1/2$ and apply the estimate (\ref{eq:4.5}) with $k = 2$
and (4.3) for $\delta$. This gives
$$|I_3| \leq \epsilon_n \lambda_{\xi}^n A_{\delta} e^{c_0|\xi|}\frac{\rho^n}{\epsilon_n^2} 
\int_{|u| > b} |u|^{\delta - 2} du = B \epsilon_n \lambda_{\xi}^n \Bigl(\frac{\rho^n}{\epsilon_n^2}\Bigr).$$
According to the condition (1.5), we can arrange for $n \geq n_0$ and $0 \leq \alpha_0 < -\frac{\log \rho}{2}$ the inequality
$$-\frac{\log n}{n} + \frac{2 \log \epsilon_n}{n} = \frac{\log n}{n} + \frac{2 \log \delta_n}{n} \geq \log \rho ,$$
which leads to
$$\frac{\rho^n}{\epsilon_n^2} \leq \frac{1}{n},\: n \geq n_0\;.$$
Thus, we conclude that 
$$I_3 = {\mathcal O} \Bigl(\frac{\epsilon_n \lambda_{\xi}^n}{n}\Bigr).$$
Consequently,
$$\rho(n) = I_1 + I_2 + I_3 =  \frac{\hat{\chi}(0)}{\sigma} \frac{\epsilon_n}{\sqrt{ 2 \pi n}} \lambda_{\xi}^n
\Bigl( 1 + {\mathcal O}(1/\sqrt{n})\Bigr)\;,$$
and taking into account (1.5), we deduce
$$\lim_{n \to \infty} \frac{1}{n} \log \rho(n) = \log \lambda_{\xi} + \lim_{n \to \infty} \frac{\log \epsilon_n}{n}  
= \Pr(\Phi + \xi_p \Psi_p) - \alpha_0 = - J (p) - \alpha_0\;.$$
This proves Proposition 1.\\

\begin{rem} By using the argument in \cite{kn:PoS1}, \cite{kn:PeS2}, under the assumptions of Theorem 1 we can show that as $n \to \infty$ we have
\begin{equation}
 m \left( \left\{ x\in R : \frac{G^{\tau^n(x)}(x)}{n}
\in (p - \delta_n , p + \delta_n)\right\}\right) \sim  \frac{2}{\sqrt{ 2 \pi n} \sigma} \exp(- n (J(p) + \alpha_0)),\:\forall p \in \Int(\II_0). 
\end{equation}
\end{rem}

\def\W{{\mathcal W}}
\def\w{{\sf w}}
\def\sAA{\Sigma_{\aa}^+}
\def\sA{\Sigma_{\aa}}
\def\ff{{\mathcal F}}

\def\hm{h^{(m)}}
\def\hmo{h^{(m-1)}}
\def\ho{h^{(0)}}
\def\oh{\overline{h}}
\def\zm{z^{(m)}}
\def\zo{z^{(0)}}
\def\zmo{z^{(m-1)}}

\section{Proof of Theorem 2}
\setcounter{equation}{0}

\subsection{Temporal functions}
\setcounter{equation}{0}

Let $\varphi_t : M \longrightarrow M$ be a $C^2$ Axiom A flow on a
$C^2$ complete (not necessarily compact) Riemannian manifold $M$ and let $\mt$ be a basic set for $\varphi_t$.
Let $\rr = \{R_i\}_{i=1}^k$ be a Markov family for $\varphi_t$ over $\mt$ as in Sect. 2.
We will assume that the size $\chi > 0$ of the Markov family $\rr$ is less than $\ep_1$, 
so that $[x,y]$ is well-defined for $x, y\in R_i$ for any $i$. 

Next, assume that $G > 0$ on $\mt$ and set $A = \min_{x\in \mt} G(x) > 0$.

Define the {\it temporal function}\footnote{It might be called the {\it temporal $\Psi$-function}, 
since it  relates in a natural way to $\Psi$.} $\Delta_\Psi$ by 
\be
\Delta_{\Psi}(x,y) = \int_0^{\Delta(x,y)} G (\varphi_t([x,y]))\, dt
\ee
for $x,y\in \mt$, $d(x,y) < \ep_1$. Just like $\Delta$, this function is constant on stable leaves 
with respect to the first variable and constant on unstable leaves with respect to the second. That is,
$\Delta_{\Psi}(x',y') = \Delta_{\Psi}(x,y)$ for any $x' \in W^s_{\ep_0}(x)\cap \mt$, $y'\in W^u_{\ep_0}(y)\cap \mt$
with $d(x',y') < \ep_1$. 

An important relationship between between $\Psi$ and $\Delta_{\Psi}$  is the following.

\begin{lem} Assume that $G \in C^\alpha(U)$ for some $\alpha > 0$.
There exists a  constant $C_1 \geq 1$ (depending on $\alpha$) such that
for any $i = 1, \ldots,k$, any $x,y \in \hU_i$ and any integer $m \geq 1$, if
$x,y$ belong to the same cylinder of length $m$ (i.e. $\sigma^j(x), \sigma^j(y)$ belong to the same 
$U_{i_j}$ for all $j = 0,1, \ldots, m-1$), then
\begin{equation}
\left|\Psi^m(x) - \Psi^m(y) + \Delta_{\Psi}(\pp^m(x), \pp^m(y))  \right| \leq
C_1\, |G|_\alpha\, (d(\pp^m(x), \pp^m(y)))^\alpha\;.
\end{equation}
\end{lem}


\noindent
{\it Proof.}  Let $x,y \in \hU_i$ and $m \geq 1$ be such that $x,y$ belong to the same cylinder of 
length $m$. 
Assume $\tau^m(x) \geq \tau^m(y)$ (the other case is similar). Since $\varphi_{\tau^m(y)}(x)\in W^u_\ep(\pp^m(y))$,
it follows that $[\pp^m(x), \pp^m(y)] = \pp^m(x)$ and $\Delta(\pp^m(x), \pp^m(y)) = -(\tau^m(x)-\tau^m(y))$.
Thus, changing the variable $t$ to $s = \tau^m(x) -t$ in the second integral in the right-hand-side below we get
\begin{eqnarray*}
\int_0^{\tau^m(x)} G  (\varphi_t(x))\, dt 
& = & \int_0^{\tau^m(y)} G (\varphi_t(x))\, dt + \int_{\tau^m(y)}^{\tau^m(x)} G (\varphi_t(x))\, dt\\
& = &  \int_0^{\tau^m(y)} G (\varphi_t(x))\, dt - \int_{\tau^m(x) -\tau^m(y)}^{0} G (\varphi_s(\pp^m(x)))\, ds\\
& = &  \int_0^{\tau^m(y)} G (\varphi_t(x))\, dt - \int_0^{\Delta(\pp^m(x), \pp^m(y))} G (\varphi_s([\pp^m(x), \pp^m(y)]))\, ds\\
& = &  \int_0^{\tau^m(y)} G (\varphi_t(x))\, dt - \Delta_{\Psi}(\pp^m(x), \pp^m(y))\;.
\end{eqnarray*}
This implies
$$|\Psi^m(x) - \Psi^m(y) +  \Delta_{\Psi}(\pp^m(x), \pp^m(y))| 
= \left| \int_0^{\tau^m(y)} \left[G(\varphi_t(x)) -  G(\varphi_t(y))\right]\, dt\right| .$$

For $D = d(\pp^m(x), \pp^m(y))$ and $0 \leq t \leq \tau^m(y)$ we have
$$|G (\varphi_t(x)) - G (\varphi_t(y)| \leq C\, |G|_\alpha\, (d(\varphi_t(x), \varphi_t(y)))^\alpha
\leq C\, |G|_\alpha\, \frac{1}{c_0\, \gamma^{\alpha(\tau^m(y) -t)}} \, D^\alpha$$
for some constants $c_0 > 0$, $\gamma > 1$ (see (2.2) for the case $\alpha = 1$), so it follows that
\begin{eqnarray*}
\left| \int_0^{\tau^m(y)} \left[ G (\varphi_t(x)) - G (\varphi_t(y))\right]\, dt\right|
& \leq & \frac{C\, |G|_\alpha\, D^\alpha}{c_0} \int_0^{\tau^m(y)} \gamma^{\alpha(t-\tau^m(y))}\, dt\\
&  =    & \frac{C\, |G|_\alpha\, D^\alpha}{c_0\, \alpha\,\log \gamma} (1- 1/\gamma^{\alpha\tau^m(y)})\\
& \leq & C_1\, |G|_\alpha\, (d(\pp^m(x), \pp^m(y)))^\alpha\;
\end{eqnarray*}
for some constant $C_1 > 0$. This proves (5.2).
\endofproof

\def\MM{{\mathcal M}}
\def\mm{{\mathcal M}}
\def\fa{f^{(a)}}
\def\fb{f^{(b)}}
\def\f0{f^{(0)}}
\def\fab{f^{(ab)}}
\def\mab{\mm_{ab}}
\def\ma{\mm_{a}}
\def\lab{\lc_{ab}}
\def\psib{\psi^{(b)}}

\def\fc{f^{(c)}}
\def\fac{f^{(ac)}}
\def\mac{\mm_{ac}}
\def\lac{\lc_{ac}}
\def\labc{\lc_{abc}}
\def\psic{\psi^{(c)}}

\def\Con{\mbox{\rm Const}}
\def\hnu{\hat{\nu}}
\def\hrho{\hat{\rho}}

\subsection{Spectral estimates for Ruelle transfer operators}

Throughout we assume that $\varphi_t$  and $\mt$ satisfy the Standing Assumptions stated in Sect. 1.

Given a Lipschitz real-valued function $f$  on $\hU$, set $\tf  = f - P\Psi$, where 
$P = P_f\in \R$ is the unique 
number such that the topological pressure $\Pr_\sigma(\tf)$ of $\tf$ with respect to $\sigma$ is 
zero (cf. e.g. \cite{kn:PP}). For $a, b\in \R$, one defines the {\it Ruelle transfer operator}
$$\lc_{\tf-(a+\i b)\Psi} : \clip (\hU) \longrightarrow \clip (\hU)$$ 
in the usual way (cf. Sect. 2 above). 

To prove Theorem 2 we apply the arguments from Sects. 3 and 5 in \cite{kn:St2} and a modification of 
the arguments in Sect. 4 of \cite{kn:St2}.
The main step is to prove the analogue of Lemma 4.3 in \cite{kn:St2} for the roof function 
$\tau$ replaced by $\Psi$. 
Below we sketch the arguments required to achieve this following the reasoning in Sect. 4 in 
\cite{kn:St2} with some modifications.

We will need the following lemma the first part of which is Lemma 4.1 in \cite{kn:St2}. 


\begin{lem} Let $i = 1, \ldots, k$, $\tz\in U_i$ and $\delta > 0$. Then there exists 
$\ep' > 0$ such that 
for any $y_1 \in W^s_{\ep}(\tz)\cap \mt$ and any $y_2 \in W^s_{\ep}(\tz)\cap \mt$ sufficiently 
close to $y_1$, we have:

\ms

{\rm (a)} {\rm (\cite{kn:St2})} $|\Delta(z',\pi_{y_2}(z)) - \Delta(z', \pi_{y_1}(z))| < \delta\, d (z,z')\;$
{\it for all $z,z' \in \mt \cap W^u_{\ep'}(\tz)$.}

\ms

{\rm (b)} $|\Delta_{\Psi}(z',\pi_{y_2}(z)) - \Delta_{\Psi}(z', \pi_{y_1}(z))| \leq
\|G\|_0\, |\Delta(z',\pi_{y_2}(z)) - \Delta(z', \pi_{y_1}(z))|$.
\end{lem}


\noindent
{\it Proof of part} (b).  Let $y_1, y_2 \in W^s_{\ep}(\tz)\cap \mt$ and $z,z' \in \mt \cap W^u_{\ep}(\tz)$.
By definition, $\pi_{y_j}(z) = [z,y_j] \in W^s_\ep(z)$ and $[z', \pi_{y_j}(z)] \in W^s_\ep(z')$. Since 
$G$ is constant on any stable leaf in any rectangle $R_i$, it follows that 
$G([z', \pi_{y_j}(z)]) = G(z')$ and moreover $G(\varphi_t([z', \pi_{y_j}(z)])) = G(\varphi_t(z'))$
for all $t \in [0,\tau(z')]$. Thus,
\begin{eqnarray*}
\Delta_{\Psi}(z',\pi_{y_2}(z)) - \Delta_{\Psi}(z', \pi_{y_1}(z))
& = &  \int_0^{\Delta(z',\pi_{y_2}(z))} G (\varphi_t([z',\pi_{y_2}(z)]))\, dt\\
&    & -  \int_0^{\Delta(z',\pi_{y_1}(z))} G (\varphi_t([z',\pi_{y_1}(z)]))\, dt\\
& = &  \int_0^{\Delta(z',\pi_{y_2}(z))} G (\varphi_t(z'))\, dt -  
\int_0^{\Delta(z',\pi_{y_1}(z))} G (\varphi_t(z'))\, dt\\
& = & \pm \int_0^{\Delta(z',\pi_{y_2}(z)) - \Delta(z',\pi_{y_1}(z)) } G (\varphi_t(z'))\, dt\;.
\end{eqnarray*} 
Hence
$$\left| \Delta_{\Psi}(z',\pi_{y_2}(z)) - \Delta_{\Psi}(z', \pi_{y_1}(z))\right|
\leq \left| \Delta_{\Psi}(z',\pi_{y_2}(z)) - \Delta_{\Psi}(z', \pi_{y_1}(z))\right|\, \|G\|_0\;.$$
This proves the lemma.
\endofproof

\bs

Following Sect. 4 in \cite{kn:St2}, {\bf fix  an arbitrary point $z_0 \in \mt$ and constants 
$\ep_0 > 0$ and $\theta_0 \in (0,1)$ 
 with the properties described  in} (LNIC). Assume that  
 $z_0 \in \Int_\mt(U_1)$, $U_1 \subset \mt \cap W^u_{\ep_0}(z_0)$ and 
$S_1 \subset \mt \cap W^s_{\ep_0}(z_0)$. Fix an arbitrary constant $\theta_1$ such that
$$0 <  \theta_0  < \theta_1 < 1 \;. $$

Next, fix an arbitrary orthonormal basis $e_1, \ldots, e_{n}$ in $E^u (z_0)$ and a $C^1$
parametrization $r(s) = \exp^u_{z_0}(s)$, $s\in V'_0$, of a small neighbourhood $W_0$ of $z_0$ in 
$W^u_{\ep_0} (z_0)$ such that $V'_0$ is a convex compact neighbourhood of $0$ in 
$\R^{n} \approx \mbox{\rm span}(e_1, \ldots,e_n) = E^u(z_0)$. Then $r(0) = z_0$ and
$\frac{\partial}{\partial s_i} r(s)_{|s=0} = e_i$ for all $i = 1, \ldots,n$.  Set 
$U'_0 = W_0\cap \mt $.
Shrinking $W_0$ (and therefore $V'_0$ as well)
if necessary, we may assume that $\overline{U'_0} \subset \Int_\mt (U_1)$ and
$\left|\left\la \frac{\partial r}{\partial s_i} (s) ,  \frac{\partial r}{\partial s_j} (s) \right\ra
- \delta_{ij}\right| $ is uniformly small for all $i, j = 1, \ldots, n$ and $s\in V'_0$, so that
\be
\frac{1}{2} \la \xi , \eta \ra \leq \la \; d r(s)\cdot \xi \; , \;  
d r(s)\cdot \eta\;\ra \leq 2\, \la \xi ,\eta \ra \quad, \quad
\xi , \eta \in E^u(z_0) \:,\: s\in V'_0\, ,
\ee
and
\be
\frac{1}{2}\, \|s-s'\| \leq  d ( r(s) ,  r(s')) \leq 2\, \|s-s'\|\quad, \quad s,s'\in V'_0\;.
\ee

\bs

\noindent
{\bf Definitions} (\cite{kn:St2}): (a)
For a cylinder $\cc \subset U'_0$ and a unit vector $\xi \in E^u(z_0)$
we will say that a {\it separation by a $\xi$-plane occurs} in $\cc$ if there exist $u,v\in \cc$ with 
$d(u,v) \geq \frac{1}{2}\, \diam(\cc)$ such that
$ \left\la \frac{r^{-1}(v) - r^{-1}(u)}{\| r^{-1}(v) - r^{-1}(u)\|}\;,\; \xi \right\ra  \geq \theta_1\;.$

Let  $\ss_\xi$ be the {\it family of all cylinders} $\cc$ contained in $U'_0$ such that a separation by an $\xi$-plane 
occurs in $\cc$.

\ms

(b) Given an open subset $V$ of $U'_0$  which is a finite union of open cylinders and  $\delta > 0$, let
$\cc_1, \ldots, \cc_p$ ($p = p(\delta)\geq 1$) be the family of maximal closed cylinders in $\oV$ with
$\diam(\cc_j) \leq \delta$. For any unit vector $\xi \in E^u(z_0)$ set
$M_{\xi}^{(\delta)}(V) = \cup \{ \cc_j : \cc_j \in \ss_{\xi} \:, \: 1\leq j \leq p\}\;.$

\ms

In what follows we will construct, amongst other things, a sequence of unit vectors
$\xi_1, \xi_2, \ldots, \xi_{j_0}\in E^u(z_0)$. For each $\ell = 1, \ldots,j_0$  set 
$B_\ell = \{ \eta \in \sn : \la \eta , \xi_\ell\ra  \geq \theta_0\}\;.$
For $t \in \R$ and $s\in E^u(z_0)$ set
$I_{\eta ,t} g(s) = \frac{g(s+t\, \eta) - g(s)}{t}$, $ t \neq 0\;$
({\it increment} of $g$ in the direction of $\eta$).

As in \cite{kn:St2}, where we dealt with the case $\Psi = \tau$, the main aim 
will be to prove the following:


\begin{lem}
There exist integers $1 \leq n_1 \leq N_0$ and $\ell_0 \geq 1$,
a sequence of unit vectors $\eta_1, \eta_2, \ldots, \eta_{\ell_0}\in E^u(z_0)$
and a non-empty open subset $U_0$ of $U'_0$ which is a finite union of open cylinders of 
length $n_1$ such that setting $\uu = \sigma^{n_1} (U_0)$ we have:

{\rm (a)} {\it For any integer $N\geq N_0$ there exist Lipschitz maps $\vl_1, \vl_2 : U \longrightarrow U$ 
($\ell = 1,\ldots, \ell_0$)  such that $\sigma^N(\vl_i(x)) = x$  for all $x\in \uu$ and $\vl_i (\uu)$ is
a finite union of open cylinders of length $N$ ($i=1,2$; $\ell = 1,2, \ldots,\ell_0$).} 

{\rm (b)} {\it There exists a constant $\hd > 0$ such that for all  $\ell = 1, \ldots, \ell_0$, 
$s\in r^{-1}(U_0)$, $0 < |h| \leq \hd$ and $\eta \in B_\ell$ so that  
$s+h\, \eta \in r^{-1}(U_0\cap \mt)$ we have }
$$\left[I_{\eta,h} \left(\Psi^{N}(\vl_2(\trr(\cdot ))) - \Psi^{N}(\vl_1(\trr(\cdot)))\right)\right](s)  \geq \frac{A \hd}{4}\,.$$

{\rm (c)} {\it We have
$\overline{\vl_i (U)} \bigcap \overline{v_{i'}^{(\ell')}(U)} = \e$ whenever $(i,\ell) \neq (i',\ell')$.}

{\rm (d)} {\it For  any open cylinder $V$ in  $U_0$ there exists a constant 
$\delta' = \delta'(V) > 0$  such that
$$V \subset M_{\eta_1}^{(\delta)}(V) \cup M_{\eta_2}^{(\delta)}(V) \cup \ldots 
\cup M_{\eta_{\ell_0}}^{(\delta)}(V)$$ 
for all} $\delta \in (0,\delta']\;.$
\end{lem}


Clearly, if $U_0$ and $\uu$ are as in the lemma, then we must have 
$\overline{\uu} = U$.

The following lemma is proved in Sect. 4 in \cite{kn:St2}.


\def\Ulo{U^{(\ell_0)}}

\begin{lem} {\rm (Lemma 4.4 in \cite{kn:St2})}
There exist an integer $\ell_0 \geq 1$, 
open cylinders $\Ulo_0 \subset \ldots \subset \Uo_0$ contained in $U'_0$, and  for each 
$\ell = 1, \ldots, \ell_0$, an integer $m_\ell \geq 1$  such that the following hold:

{\rm (i)} {\it For each 
$i = 1,2$,  there exists a contracting map
$\wl_i : \Ul_0 \longrightarrow U_1$ such that $\sigma^{m_\ell}(\wl_i(x)) = x$ for all $x\in \Ul_0$,
$\wl_i : \Ul_0 \longrightarrow \wl_i(\Ul_0)$ is a homeomorphism, 
$\wl_i(\Ul_0)$ is an open cylinder of length $\geq m_\ell$ in $U_1$, and 
the  sets $\overline{\wl_i(\Ul_0)}$ ($\ell = 1, \ldots, \ell_0$, $i = 1,2$)
are disjoint.}

{\rm (ii)} {\it  For each $\ell = 1, \ldots, \ell_0$ there exist a number $\hd_\ell \in (0,\delta_0)$ 
and a vector $\eta_\ell\in \sn$ such that
$| \langle  \eta_\ell , \eta_{\ell'}  \rangle |  \leq  \theta_1$ whenever $\ell' \neq \ell\;,$ and
$$\inf_{0 < |h| \leq \hd_\ell} \left|\left[I_{\eta ,h}\left( \tau^{m_\ell}(\wl_2(r(\cdot))) 
- \tau^{m_\ell}(\wl_1(r(\cdot))) \right)\right](s)\right| \geq \hd_\ell \quad, \quad r(s)\in \Ul_0\:, \: \eta \in B_\ell \:,$$
for all $\ell = 1, \ldots, \ell_0$, where the $\inf$ is taken over $h$ with $r(s+h \eta )\in \Ul_0$. }

{\rm (iii)} {\it For  any open cylinder $V$ in  $\Ulo_0$ there exists a constant 
$\delta' = \delta'(V)\in (0,\delta_0)$  such that
$V \subset M_{\eta_1}^{(\delta)}(V) \cup M_{\eta_2}^{(\delta)}(V) \cup \ldots 
\cup M_{\eta_{\ell_0}}^{(\delta)}(V)$ for all} $\delta \in (0,\delta']\;.$
\end{lem}

In what follows we use the objects constructed in Lemma 5.
Set 
$\di \hd = \min_{1\leq \ell\leq \ell_0} \hd_j$, 
$\di n_0 = \max_{1\leq \ell\leq \ell_0} m_\ell$, 
and fix an arbitrary point $\hz_0 \in U_0^{(\ell_0)}\cap \hU$. 
Set 
$$L = \Lip(G) \quad ,   \quad D = \diam(\mt)\,,$$
and recall that $A = \min_{x\in \mt} G(x) > 0$.
Without loss of generality we may assume that $D \geq 1$. We will also assume that
\be
\mu_0 \leq \frac{c_0\, \hd}{128\, C_0\, C_1\, D}\;.
\ee

Next, essentially repeating some arguments from the proof of Lemma 5 we get
an analogue of property (ii) there with $\tau$ replaced by $\Psi$.

\begin{lem}
Assume that {\rm (5.5)} holds and $L \leq \mu_0\, A$. Then for every 
$\ell = 1, \ldots, j_0$ we have
$$\inf_{0 < |h| \leq \hd_\ell} | I_{\eta ,h} [\Psi^{m_\ell}(\wl_2(r(s))) - \Psi^{m_\ell}(\wl_1(r(s)))] | \geq \frac{A\, \hd_\ell}{4} 
\quad,\quad r(s)\in \Ul_0\:, \: \eta \in B_\ell \:,$$
where the $\inf$ is taken over $h$ with $r(s+h \eta )\in \Ul_0$. 
\end{lem}

\noindent
{\it Sketch of proof of Lemma} 6.  
Clearly $\|G\|_0 \leq A+LD$, so assuming $L \leq \mu_0\,A$, we have $\|G\|_0\leq A( 1+ \mu_0\, D)$. 

We will follow the inductive steps in the proof of Lemma 4.4 in \cite{kn:St2}. In fact we will do in
detail the case $\ell = 1$; in the general case the argument is very similar.

\ms

\noindent
{\bf Step 1.} As in \cite{kn:St2}, consider a unit vector $\xi_1 \in E^u(z_0)$ tangent to $\mt$ at $z_0$, and
then using the condition (LNIC) and the choice of $z_0$ that there exist $\tz = r(\ts) \in U'_0$, 
$y_1', y_2' \in W^s_{R_1}(\tz)$ with $y'_1 \neq y'_2$, $\delta'_1 > 0$, 
$\ep'_1 > 0$, and an open cylinder $\Uo_0$ contained in $B_U(\tz,\ep'_1) \subset U'_0$ 
with $\tz \in \Uo_0$ such that
$$|\Delta (r(s + h\, \eta), \pi_{y'_1}(r(s))) - \Delta (r(s + h\, \eta), \pi_{y''_1}(r(s)))| 
\geq \hd_1\, |h| \;$$
whenever $r(s) \in \Uo_0$, $ \eta \in B_1$, $| h| \leq \hd_1$ and $r(s+h \eta) \in \Uo_0$,
where $\hd_1 = \min\{ \delta'_1/2, \ep'_1\}$. 
The above and Lemma 3(b) imply 
\begin{equation} 
|\Delta_\Psi (r(s + h\, \eta), \pi_{y'_1}(r(s))) - \Delta_\Psi (r(s + h\, \eta), \pi_{y''_1}(r(s)))|  \geq A \hd_1\, |h| .
\end{equation}

Now  for $z = r(s) \in \Uo_0$, $\eta\in B_1$ and $h\in \R$ with $r(s + h\, \eta) \in \Uo_0$ using (5.2)  we get
\begin{eqnarray}
&          &\left| [\Psi^{m_1}(\wo_2(r(s+ h\, \eta))) - \Psi^{m_1}(\wo_1(r(s+ h\, \eta)))]
  - [\Psi^{m_1}(\wo_2(r(s))) - \Psi^{m_1}(\wo_1(r(s)))]\right|\nonumber\\
& =       & \left| [\Psi^{m_1}(\wo_1(r(s))) - \Psi^{m_1}(\wo_1(r(s+ h\, \eta)))]
- [\Psi^{m_1}(\wo_2(r(s)))  - \Psi^{m_1}(\wo_2(r(s+ h\, \eta)))] \right|\nonumber\\
& \geq  & \left| \Delta_\Psi ( \pp^{m_1}(\wo_1(r(s+ h\, \eta))), \pp^{m_1}(\wo_1(r(s)) )
  - \Delta_\Psi( \pp^{m_1}(\wo_2(r(s + h\, \eta))) , \pp^{m_1}(\wo_2(r(s)) )\right|\nonumber\\
&          & - C_1 \, L\, \sum_{i=1}^2d( \pp^{m_1}(\wo_i(r(s+ h\, a))), \pp^{m_1}(\wo_i(r(s)) )\nonumber\\
&     =   & \left| \Delta_\Psi ( \pi_{y'_1}(r(s+ h\, \eta)), \pi_{y'_1}(r(s)) ) 
- \Delta_\Psi (\pi_{y''_1}(r(s + h\, \eta) ),  \pi_{y''_1}(r(s) ))\right| \nonumber\\
&          & - C_1 \, L\, \sum_{i=1}^2d( \pp^{m_1}(\wo_i(r(s+ h\, \eta))), \pp^{m_1}(\wo_i(r(s)) )\,.
\end{eqnarray}
Since $\pp^{m_1}(\wo_1(x)) = \pi_{y'_1}(x)$ and $\pp^{m_1}(\wo_2(x)) = \pi_{y''_1}(x)$, 
using (2.1) and (5.4) we get 
\begin{eqnarray*}
d( \pp^{m_1}(\wo_1(r(s+ h\, \eta))), \pp^{m_1}(\wo_1(r(s)) )
& =    & d(\pi_{y'_1}(r(s+ h\, \eta)), \pi_{y'_1}(r(s)) )\\
& \leq & C_0 \, d(r(s+ h\, \eta)), r(s) ) \leq 2C_0 \, |h|\;.
\end{eqnarray*}
Sinilarly, $d( \pp^{m_1}(\wo_2(r(s+ h\, \eta))), \pp^{m_1}(\wo_2(r(s)) ) \leq 2C_0 \, |h|$.
This, (5.7), (5.6) and (5.5)  imply
\begin{eqnarray*}
&          & \left| I_{\eta,h} [\Psi^{m_1}(\wo_2(r(s))) - \Psi^{m_1}(\wo_1(r(s)))] \right| \\
& \geq  & \left| \Delta_\Psi ( r(s+ h\, \eta), \pi_{y'_1}(r(s)) ) - 
\Delta_\Psi (r(s + h\, \eta) ,  \pi_{y''_1}(r(s) ))\right| 
           - 4 C_0 C_1 L\, |h| \\
& \geq  & A\,  \hd_1\, |h| - 4 \mu_0 A \, C_0 \, C_1\, |h|  \geq \frac{A\, |h|\, \hd_1}{2}
\end{eqnarray*}
whenever $r(s)\in \Uo_0$, $\eta \in B_1$, $0 < |h| \leq \ep'_1$ and $r(s+h\, \eta) \in \Uo_0$.

This proves the statement for $\ell = 1$. As we mentioned above, in the general case the argument
is the same.
\endofproof

\bs

Next, as in \cite{kn:St2} we {\bf fix  $n_1 > 0$ and an open neighbourhood $U_0$ of $\hz_0$ in 
$U_0^{(j_0)}$} such that
$U_0$ is a finite union of open cylinders of length $n_1$, $\uu = \sigma^{n_1}(U_0) = \Int_{\mt}(U )$
and $\sigma^{n_1} : U_0 \longrightarrow  \uu$ is a homeomorphism.
The inverse homeomorphism $\psi : \uu \longrightarrow U_0$ is Lipschitz, so it has a Lipschitz extension
\begin{equation}
\psi : U \longrightarrow \overline{U_0}\:\: \mbox{\rm  such that }\:\: 
\sigma^{n_1}(\psi(x)) = x \:\:, \:\:x\in \uu\;.
\end{equation}
Then $\trr(s) = \sigma^{n_1}(r(s))$,  $s\in V_0\;,$ where $V_0 = r^{-1}(U_0) \subset V'_0$,
gives a Lipschitz parametrization of $\uu$ with $\psi(\trr(s)) = r(s)$ for all $s \in V_0$.
Finally, set
\be
\Vj_i = \wj_i(U_0)  \subset \tVj_i\quad, \quad i = 1,2\; ; \: j = 1, \ldots,j_0\;.
\ee
It follows from the choice of $U_0$, the properties of $\wj_i$ (see (i) in Lemma 5 and Proposition 3.1
in \cite{kn:St2} that $\Vj_i$ is a finite union of open cylinders of lengths $n_1 + m_j$.

The following  lemma is proved by using arguments from  \cite{kn:D} (see also \cite{kn:St2}). 
Since the proof of the analogous lemma in \cite{kn:D} and \cite{kn:St2} 
is very sketchy (or absent), we provide some details here.

\ms

\begin{lem} 
Let $\delta '' > 0$ and assume that $\mu_0 \leq \frac{\delta''}{4C_0C_1}$. Then there exists an integer 
$n_2 > 0$ such that for any $m  \geq n_0+ n_2$, 
any $j = 1, \ldots, j_0$ and $i = 1,2$ there exist contracting maps  
$\tvj_i : \Vj_i  \longrightarrow U $  with $\sigma^{m-m_j}(\tvj_i(w)) = w$ for all $w\in \Vj_i $ such that
\be
\Lip (\Psi^{m-m_j} \circ \tvj_i )\leq A\, \delta'' \:\:\: \mbox{\rm   on }\:\:\:  \Vj_i \;,
\ee
$\tvj_i (\Vj_i )$ is a finite union of open cylinders of length $n_1 + m$ and
$\overline{\tvj_i (\Vj_i )} \bigcap \overline{\tilde{v}_{i'}^{(j')}(V_{i'}^{(j')})} = \e$ whenever
$(i,j) \neq (i',j')$.
\end{lem}

\ms

\noindent
{\it Proof of Lemma} 7. It follows from Lemma 3 and $\|\Phi\|_0 \leq A+LD \leq A(1+\mu_0 D)$ 
that there exists $\delta > 0$ such that for any $y\in R_1$ with $d (y,\piU(y)) < \delta$  we have 
$$|\Delta_\Psi(x',\pi_y(x))| \leq  \frac{A\,\delta''}{2}\, d (x,x')\;$$ 
for all  $x, x'\in U_1$. 

Fix  arbitrary non-empty open subset (e.g. cylinders) $Q_1, \ldots, Q_{j_0}$  in $U_1$ 
such that $\overline{Q_j} \cap \overline{Q_{j'}} = \e$ whenever $j \neq j'$.
Now take $n_2$ so large that for any $m \geq n_2$ and any $j = 1, \ldots,j_0$,
$\pp^m(Q_j)$ contains a whole leaf $W^u_{R_1}(z)$ for some $z\in R_1$ such that 
$d (y,\piU(y)) < \delta$ for all $y \in W^u_{R_1}(z)$.
This is possible, since $\pp^m(Q_j)$ fills in $R_1$ densily as $m \to \infty$. 

Let $m \geq n_2 + n_0$. Given $j = 1, \ldots,j_0$, we have $m-m_j \geq n_2$,
so we can choose $\zeta_j\in R_1$ such that $W^u_{R_1}(\zeta_j) \subset \pp^{m-m_j}(Q_j)$ and
$d (y,\piU(y)) < \delta$ for all $y \in W^u_{R_1}(\zeta_j)$. There exists
a  (closed) cylinder $L_j$ in $Q_j \subset U_1$ such that the Lipschitz extension of the restriction  of $\pp^{m-m_j}$ 
to $\Lo_j$ defines a Lipschitz homeomorphism between $L_j$ and $W^u_{R_1}(z_j)$, so in particular
$\sigma^{m-m_j}$ can be extended from $\Lo_j$ to a Lipschitz homeomorphism  $\sigma^{m-m_j} : L_j \longrightarrow U_1$. 
Then for any $w\in U_1$ there exists a unique element $\tvj_i(w)$ of $L_j$ with
$\piU (\pp^{m-m_j}(\tvj_i(w))) = w$, i.e. $\sigma^{m-m_j}(\tvj_i(w)) = w$. It is clear from the
definition that $\tvj_i : U_1 \longrightarrow L_j \subset U_1$ is contracting (assuming that  $n_2$ is sufficiently large). 
For $w,w'\in U_1$, as in the the proof of (5.7) we have
\begin{eqnarray*}
&      & |\Psi^{m-m_j}(\tvj_i(w')) - \Psi^{m-m_j}(\tvj_i(w)) |\\
& \leq & |\Delta_\Psi (\pp^{m-m_j}(\tvj_i(w)), \pp^{m-m_j}(\tvj_i(w')))| 
+ C_1L\, d((\pp^{m-m_j}(\tvj_i(w)), \pp^{m-m_j}(\tvj_i(w'))\\
& \leq & |\Delta_\Psi( w, \pi_{z_j}(w') )|+ 2 C_0C_1L\, d(w,w') \leq 
A\, \left[ \delta''/2 + 2C_0C_1\, \mu_0\right]\, d (w,w') \leq A\, \delta''\, d(w,w')\;.
\end{eqnarray*}
In particular, $\Lip( \Psi_{m-m_j} \circ \tvj_i)  \leq A\, \delta''$. To conclude the proof, 
consider the restriction of the 
map $\tvj_i$ to $\Vj_i$; then by Proposition 3.1 in \cite{kn:St2}, $\tvj_i (\Vj_i )$ is a finite 
union of open cylinders of lengths 
$n_1+m_j + (m-m_j) = n_1 + m$. 

It remains to show that the cylinders $\tvj_i (\Vj_i )$ are disjoint. Assume that 
$\overline{\tvj_i (\Vj_i )}$
and $\overline{\tilde{v}_{i'}^{(j')}(V_{i'}^{(j')})}$ have a common point for some 
$i,i' = 1,2$ and $j, j' = 1,2,\ldots,j_0$.
Since $\tvj_i (\Vj_i ) \subset L_j$ and $\tilde{v}_{i'}^{(j')}(V_{i'}^{(j')}) \subset L_{j'}$, it follows
that $L_j \cap L_{j'} \neq \e$, so $Q_j \cap Q_{j'} \neq \e$. This implies
$j = j'$ and therefore $m_j = m_{j'}$.  It now follows that
$\Vj_i = \sigma^{m- m_j}(\tvj_i (\Vj_i ))$ and $V_{i'}^{(j)} 
= \sigma^{m-m_j}(\tilde{v}_{i'}^{(j)}(V_{i'}^{(j)}) )$
have a common point, which implies $i = i'$.
$\endofproof$

\bigskip

Set $\di \delta'' = \frac{c_0\, \hd}{32}$ and notice that (5.5) implies 
$\di \mu_0 \leq \frac{\delta''}{4C_0C_1}$, so the assumptions of Lemma 7 are satisfied.
{\bf Fix $n_2 = n_2(\delta'') > 0$ with the properties listed  in Lemma 7}, and  denote
$N_0 = n_0 + n_1 + n_2\;.$


\ms

\noindent
{\it Proof of Lemma} 4. Parts (a), (c) and (d) follow from Lemma 4.2 in \cite{kn:St2} and its proof.
It remains to prove (b). We will follow the proof of Lemma 4.2 in \cite{kn:St2}.

Let $N \geq N_0$.  Then $m = N - n_1 \geq n_0 + n_2$, 
so  by Lemma 6 for any $j = 1, \ldots,j_0$ and any  $i = 1,2$ there exists a  contracting homeomorphism
\begin{equation}
\tvj_i : \Vj_i \longrightarrow \tvj_i (\Vj_i ) \subset U\:\: \mbox{\rm  with } \:\:  \sigma^{N- m_j - n_1}(\tvj_i(w)) = w \:\: , \:\: w\in \Vj_i \; ,
\end{equation}
such that (5.11) holds with $m = N - n_1$ and $\delta''$ as above. Moreover, we can choose
the maps $\tvj_i$ so that $\overline{\tvj_i (\Vj_i )} \bigcap  \overline{\tilde{v}_{i'}^{(j')}(V_{i'}^{(j')})} = \e$ whenever
 $(i,j) \neq (i',j')$. Now define Lipschitz maps
\begin{equation}
\vj_i : U  \longrightarrow U \quad \mbox{\rm such that} \quad \vj_i(x) = 
\tvj_i(\wj_i(\psi(x)))\:, \: x\in \uu\;.
\end{equation}
It follows immediately from the above that $\vj_i (U) \cap v_{i'}^{(j')}(U) = \e$ 
whenever $(i,j) \neq (i',j')$, while Proposition 3.1 in \cite{kn:St2} shows that each $\vj_i (\uu)$ is
a finite union of  open cylinders of length $N$.

Moreover, for any $x\in \uu$, according to  part (i) in Lemma 5, (5.12)  and (5.11), we have
$$\sigma^{N-n_1}(\vj_i(x)) = \sigma^{m_j}( \sigma^{N-m_j-n_1}(\vj_i(x)))  = \sigma^{m_j} (\wj_i(\psi(x))) = \psi(x)\;,$$
which is the same for all $j$ and  $i$. Consequently,
$\sigma^p(\vj_1(x)) = \sigma^p(\vj_2(x))$ for all $p \geq N-n_1$ and $x\in \uu $. Thus,
$$\Psi^N(\vj_2(x)) - \Psi^N(\vj_1(x)) = \Psi^{N-n_1}(\vj_2(x))- \Psi^{N-n_1}(\vj_1(x))$$ 
for $x \in \uu$, and given $\eta \in B_j$ and $h > 0$, we have
\begin{eqnarray*}
&    & I_{\eta ,h}\, \left[\Psi^N(\vj_2(\trr(s))) - \Psi^N(\vj_1(\trr(s)))\right]\\
& = &   I_{\eta ,h} \, \left[\Psi^{N-n_1}(\tvj_2(\wj_2(\psi(\trr(s)))))  - \Psi^{N-n_1}(\tvj_1(\wj_1(\psi(\trr(s))))) \right] \nonumber\\
& = &  I_{\eta,h}\, \left[\Psi^{N-n_1}(\tvj_2(\wj_2(r(s))))  -  \Psi^{N-n_1}(\tvj_1(\wj_1(r(s)))) \right] \;. \nonumber
\end{eqnarray*}
Since 
$$\Psi^{N-n_1}(z) = [\Psi(z) + \ldots + \Psi(\sigma^{N- m_j -n_1-1}(z))] + [\Psi(y) + \ldots + \Psi(\sigma^{m_j - 1}(y))],$$
where $y = \sigma^{N - m_j - n_1}(z)$, using (5.11) we get
\begin{eqnarray*}
\Psi^{N-n_1}(\tvj_i(\wj_i(x)))
& = & \Psi^{N- m_j - n_1}(\tvj_i(\wj_i(x))) +  \Psi^{m_j} (\sigma^{N- m_j -n_1}(\tvj_i(\wj_i(x)))\\
& = & \Psi^{N - m_j - n_1}(\tvj_i(\wj_i(x))) + \Psi^{m_j} (\wj_i(x)).
\end{eqnarray*}
Therefore, 
\begin{eqnarray*}
&   &  I_{\eta ,h}\, \left[\Psi^N(\vj_2(\trr(s))) - \Psi^N(\vj_1(\trr(s)))\right]\\
& = &   I_{\eta ,h}\, \left[\Psi^{N- m_j -n_1}(\tvj_2(\wj_2(r(s)))) -  \Psi^{N-m_j-n_1}(\tvj_1(\wj_1(r(s))))\right]\\
&    & \:\: \: +  I_{\eta,h}\, \left[\Psi^{m_j}(\wj_2(r(s))) - \Psi^{m_j}(\wj_1(r(s)))\right]\;.
\end{eqnarray*}
Let $\eta \in B_j$, and let $s$ and $s+ h\, \eta$ (for some $h$ with $0 < |h| \leq \hd$) belong to  $V_0$.
Since $\Lip ( \wj_i(x)) \leq 1/(c_0 \gamma^{m_j}) < 1/c_0$, using (5.4) and (5.10) with  $m = N-n_1$, it follows that
\begin{eqnarray*}
\left|I_{\eta ,h}\, \left[\Psi^{N- m_j -n_1}(\tvj_i(\wj_i(r(s))))\right]\,\right|
& \leq & \Lip\left[\Psi^{N-m_j-n_1}(\tvj_i)\right] \cdot \Lip (\wj_i)\cdot  \Lip( r )\\
& \leq & A\, \delta''\, \frac{ 2}{c_0} \leq \frac{A\, \hd}{8}\;.
\end{eqnarray*}
Combining this with the above and Lemma 6, gives
$$ \left| I_{\eta ,h}\, \left[\Psi^N(\vj_2(\tr(s))) - \Psi^N(\vj_1(\tr(s)))\right] \, \right|\geq 
\frac{A\,\hd}{2} - 2 \frac{A\, \hd}{8} = \frac{A\, \hd}{4}\;.$$
This proves the assertion. $\endofproof$

\bs

\noindent
{\it Proof of Theorem 2.} For Lipschitz fucntions $f$ this is just a repetition of Sect. 5 in
\cite{kn:St2} without any changes. For H\"older continuous $f$ one just needs to combine this with the approximation procedure
in \cite{kn:D} (see also Sect. 3 in \cite{kn:St1}). Since $\Psi$ is Lipschitz, the approximation procedure can be carried out
in the same way as in \cite{kn:D}.
\endofproof

\bs

\end{document}